\documentclass[a4paper,11pt]{article}
    
    \usepackage[utf8]{inputenc}
    
    \usepackage{hyperref}
    \usepackage{boxedminipage}
    
    \usepackage{amsfonts}
    
    \usepackage{amsmath} 
    
    \usepackage{amssymb}
    
    \usepackage{graphicx}
    
    \usepackage{amsthm}
    
    \usepackage{t1enc}
    
    \usepackage{subfig}
    
    \usepackage{cancel}
    
     \usepackage{textcmds}
    
     \usepackage{mathabx}

    \newtheorem{theorem}{Theorem}[section]
    
    \newtheorem{lemma}[theorem]{Lemma}

    \newtheorem{definition}[theorem]{Definition}

    \newtheorem{question}{Question}
    
    \newtheorem*{theorem*}{Theorem}

    \newcommand{\forceP}{\mathbb{P}}
    
    \newcommand{\forceQ}{\mathbb{Q}}
    
    \newcommand{\forceR}{\mathbb{R}}

    \newcommand{\ZFC}{\mathsf{ZFC}}
    
    \newcommand{\ZF}{\mathsf{ZF}}
    
    \newcommand{\ZFP}{\mathsf{ZF}^-}

    \newcommand{\BPFA}{\mathsf{BPFA}}

    \newcommand{\CH}{\mathsf{CH}}

    \newcommand{\PD}{\mathsf{PD}}
    
    \def\undertilde#1{\mathord{\vtop{\ialign{##\crcr
    
    $\hfil\displaystyle{#1}\hfil$\crcr\noalign{\kern1.5pt\nointerlineskip}
    
    $\hfil\tilde{}\hfil$\crcr\noalign{\kern1.5pt}}}}}
    
    \title{Forcing upper $\Sigma$-uniformization in the presence of lower $\Pi$-reduction or uniformization}
    
    \author{ Stefan Hoffelner\footnote{This research was funded in whole by the Austrian Science Fund (FWF) Grant-DOI 10.55776/P37228.}  }
    
\begin{document}

    \maketitle

 \begin{abstract}
 We present a method which allows the combination of forcing uniformization on the $\Pi$- and the $\Sigma$-side of the projective hierarchy to a certain extent. Using this method
we construct a universe where ${\Pi}^1_3$-reduction holds, $\Pi^1_3$-uniformization fails, yet $\Sigma^1_n$-uniformization is true for $n \ge 4$. We also construct a universe where $\Pi^1_3$-uniformization holds and for every $n \ge 4$, $\Sigma^1_n$-uniformization holds, lowering the best known upper bound for this statement from the existence of two Woodin cardinals to $Con(\ZFC)$.
\end{abstract}

\section{Introduction}

The problem of finding well-behaved choice functions for families of sets, a core concern in descriptive set theory, is a long-studied subject. Specifically, the uniformization problem, first posed by N. Lusin in 1930, asks whether one can find a choice function (called a uniformization) that has the same complexity (i.e. sits on the same projective level) as the set it is defined from.

\begin{definition}
For a set $A \subset 2^{\omega} \times 2^{\omega}$ (a set of pairs of reals), a partial function $f$ is a uniformization of $A$ if $f$ picks exactly one $y$ for each $x$ such that there is a $y'$ with $(x, y') \in A$. In other words, the graph of $f$ is a subset of $A$, and the domain of $f$ is the projection of $A$ onto its first coordinate. A projective pointclass $\Gamma \in \{ \Sigma^1_n \mid n \in \omega\}$ $ \cup \{\Pi^1_n \mid n \in \omega \}$  has the uniformization property if every set in $\Gamma$ can be uniformized by a function whose graph also belongs to $\Gamma$.
\end{definition}

Closely related to the uniformization property is a second classical regularity property, the reduction property, introduced in 1936 by K. Kuratowski.

\begin{definition}
We say that a projective pointclass $\Gamma \in \{ \Sigma^1_n \mid n \in \omega\} \cup \{\Pi^1_n \mid n \in \omega \}$ satisfies the $\Gamma$-reduction property (or just reduction) if every pair $B_0,B_1$ of $\Gamma$-subsets of the reals can be reduced by a pair of $\Gamma$-sets $R_0,R_1$, which means that $R_0 \subset B_0$, $R_1 \subset B_1$, $R_0 \cap R_1= \emptyset$ and $R_0 \cup R_1=B_0 \cup B_1$.
\end{definition}

It follows immediately from the definitions that $\Gamma$-uniformization implies ${\Gamma}$-reduction.

It is a classical result by M. Kondo that the family of $\Pi^1_1$ sets  possesses the uniformization property, which also implies the property for $\Sigma^1_2$ sets. Standard set theory ($\mathsf{ZFC}$) cannot prove uniformization for any pointclass of higher complexity.

To prove uniformization for more complex sets, one must appeal to additional set-theoretic axioms. In the constructible universe $L$, for every $n \ge 2$, $\Sigma^1_{n}$ does have the uniformization property which follows from the existence of a good $\Delta^1_2$-wellorder by an old result of Addison (see \cite{Addison}). Recall that a $\Delta^1_n$-definable wellorder $<$ of the reals is a \emph{good} $\Delta^1_n$-wellorder if $<$ is of ordertype $\omega_1$ and the relation $<_{I} \subset (\omega^{\omega})^2$ defined via
\[ x <_I y \Leftrightarrow \{(x)_n \, : \, n \in \omega\}= \{z \, : \, z < y\} \]
where $(x)_n$ is some fixed recursive partition of $x$ into $\omega$-many reals, is a $\Delta^1_n$-definable relation. Addison's observation is that a good $\Delta^1_n$-wellorder implies that the $\Sigma^1_m$-uniformization property is true for every $m \ge n$. It is easy to check that the canonical wellorder of the reals in $L$ is a good $\Delta^1_2$-wellorder so the $\Sigma^1_n$-uniformization property follows for $n \ge 2$.

On the other hand, and more importantly, if we assume the existence of large cardinals or alternatively strong forcing axioms, we obtain a very different picture.
Due to Moschovakis (see \cite{Moschovakis2} or \cite{Kechris}, Theorem 39.9), the axiom of projective determinacy ($\PD$) implies that $\Pi^1_{2n+1}$ and $\Sigma^1_{2n+2}$-sets have the uniformization property. By the Martin-Steel theorem (see \cite{MS} or \cite{Schindler}, Theorem 13.6.), the assumption of infinitely many Woodin cardinals implies $\PD$, and hence large cardinals fully settle the behaviour of the uniformization property within the projective hierarchy.

An old result of Novikov (see \cite{Schindler}, Lemma 7.25) states that $\Gamma$-reduction rules out reduction for the projective pointclass dual to $\Gamma$, denoted by $\check{\Gamma}$ so in particular the two possibilities for the behaviour of uniformization (and hence of reduction) within the projective hierarchy contradict each other.

There are set theoretic universes which combine the behaviour of uniformization under $V=L$ and under $\PD$ to some extent. Paradigmatic examples are the canonical inner models with $n$-many Woodin cardinals $M_n$. Due to Steel, these models have good $\Sigma^1_{n+2}$-definable well-orders of their reals, so in particular $\Sigma^1_{m+2}$-uniformization holds there for $m\geq n$. Yet the presence of the $n$ Woodin cardinals implies $\mathbf{\Pi}^1_n$-determinacy, so, by Moschovakis' theorem, $\Pi^1_{2m+1}$-uniformization holds for every $m$ with $2m\leq n$.

This work is a first attempt to find possibilities to combine the $\Sigma$-uniformization property and the $\Pi$-uniformization property, but using forcing techniques instead of (local forms of) projective determinacy. The articles \cite{Reduction} \cite{Uniformization} and \cite{BPFA and uniformization} deal with forcing $\Pi^1_n$-reduction, $\Pi^1_n$-uniformization for a given, fixed $n \in \omega$, and forcing global $\Sigma^1_m$-uniformization for $m \ge 2$ respectively. 
The techniques introduced there are however impossible to combine. Indeed, the methods are such that they exclude themselves from each other. Forcing uniformization (or reduction) on the $\Pi$-side of the projective hierarchy currently needs the solution of an associated fixed point problem for sets of forcings first. These fixed point problems can only be solved if we work with iterated coding forcings which satisfy some form of closure under taking products.
Forcing the uniformization property on the $\Sigma$-side according to \cite{BPFA and uniformization}, in stark contrast, needs the iteration of coding forcings which can not be closed under taking products by design (see the footnote 4 on pp. 19).

Hence it becomes a challenging problem to investigate ways to combine forcing constructions for the $\Pi$ and the $\Sigma$-uniformization property. 
Questions of this type do have a deeper underpinning. Indeed, they can be seen as first test questions for finding constructions which would yield alternative, global, or at least less local patterns for the behaviour of the uniformization property within the projective hierarchy. The latter question is still wide open \footnote{E.g. the problem of finding a universe where $\Pi^1_3$ and $\Pi^1_6$-uniformization holds is wide open and seems to be far out of reach of the current methods.} in general and the current status of available methodology suggests that it is a very hard one. This article adds further difficulties to the picture, as its methods have a degree of limited flexibility and any attempt to force, say, $\Pi^1_5$-uniformization without large cardinals must be such that it can not be combined with these methods.

The main goal of this work is the introduction of a new technique which allows the combination of forcing the $\Pi$-uniformization and methods which force the $\Sigma$-uniformization to some extent. This technique allows for a deep investigation of reduction and uniformization as displayed in the following first main theorem. Notably we do not assume anything beyond the consistency of $\mathsf{ZFC}$ for its proof.

\begin{theorem}
There is a generic extension of $L$ which satisfies
\begin{enumerate}
\item $\Pi^1_3$-reduction,
\item a failure of $\Pi^1_3$-uniformization,
\item $\Sigma^1_n$-uniformization for every $n \ge 4$ simultaneously.
\end{enumerate}
\end{theorem}

The second main theorem is another application of the ideas involved in proving the first main theorem and shows that the two Woodin cardinals of $M_2$ are not necessary to obtain the behaviour of uniformization there, instead one can get by with just assuming Con($\mathsf{ZFC})$.

\begin{theorem}
Assume the consistency of $\mathsf{ZFC}$. Then there is a universe where the $\Pi^1_3$-uniformization property is true and the $\Sigma^1_n$-uniformization property is true for all $n \ge 4$ simultaneously.
\end{theorem}

We briefly record the corresponding boldface consequences. For projective pointclasses, lightface reduction or uniformization implies the corresponding boldface property by the standard universal-set argument: one applies the lightface property to a universal set, regards the real parameter together with the first variable as one coordinate, and then takes sections. Consequently, the model of the first main theorem satisfies boldface $\boldsymbol{\Pi}^1_3$-reduction and boldface $\boldsymbol{\Sigma}^1_n$-uniformization for every $n\geq4$. By Novikov's theorem, boldface $\boldsymbol{\Sigma}^1_3$-reduction therefore fails in this model, and hence so does boldface $\boldsymbol{\Sigma}^1_3$-uniformization. This does not give the exact boldface analogue of the first main theorem: the present diagonalization yields only the failure of lightface $\Pi^1_3$-uniformization. To obtain the failure of boldface $\boldsymbol{\Pi}^1_3$-uniformization, one would have to diagonalize simultaneously against $\Pi^1_3(p)$-uniformizing graphs for every real parameter $p$ in the final model, including parameters which appear only during the generic iteration. We do not know whether the present method can be modified to achieve this, and this appears to be a technically difficult problem. In contrast, the second main theorem immediately has its full boldface analogue: its model satisfies boldface $\boldsymbol{\Pi}^1_3$-uniformization and boldface $\boldsymbol{\Sigma}^1_n$-uniformization for every $n\geq4$.

We expect that both main theorems admit corresponding higher-level analogues over the canonical inner models $M_n$, with the Woodin cardinals preserved and the projective levels shifted in the natural way. Such a lift is not a formal modification of the proofs over $L$. It requires higher-level coding predicates together with the relevant mouse-iterability, comparison, localization, and preservation arguments. Related forms of this technology are developed in \cite{Uniformization,DelfinoSigma,DelfinoPi}. Since a complete treatment would require substantial additional notions and proofs, we do not pursue these higher-level analogues in the present paper.

Put in a wider context, this work contributes to the project of a detailed investigation of separation, reduction and uniformization in the absence of large cardinals. There is a growing body of work which explores its connections, in various set theoretic contexts (see e.g. \cite{Separation without reduction}, \cite{Forcing axioms and uniformization}, \cite{MA and failure of separation}).

We end the introduction with a short description of the organization of this article. Sections~2--4 reintroduce the coding machinery and the notion of allowable forcing. Sections~5 and~6 develop the thinning-out processes, including the restart mechanism used when new forbidden information is introduced. Sections~7 and~8 construct and analyze the model for the first main theorem. Section~9 gives the modified thinning-out construction for the second main theorem and verifies its conclusions. The final mathematical section collects open questions.

    \section{Preliminaries}
    
    The techniques of this work fully rely on \cite{Reduction}, \cite{Uniformization} and \cite{Separation without reduction}. We shall give a short presentation of these now. We will indicate where the reader can find the proofs in detail.
    The coding method of our choice utilizes Suslin trees, which can be generically destroyed in an independent way of each other (see \cite{NS saturated and definable}). 
    
    \begin{definition}
    
     Let $\vec{T} = (T_{\alpha} \, : \, \alpha < \kappa)$ be a sequence of Suslin trees. We say that the sequence is an 
    
     independent family of Suslin trees if for every finite set of pairwise distinct indices $e= \{e_0, e_1,...,e_n\} \subset \kappa$ the product $T_{e_0} \times T_{e_1} \times \cdot \cdot \cdot \times T_{e_n}$ 
    
     is a Suslin tree again.
    
    \end{definition}
    
    Note that an independent sequence of Suslin trees $\vec{T} = (T_{\alpha} \, : \, \alpha < \kappa)$ has the property that if $A \subset \kappa$ and we form 
    
    $ \prod_{i \in A} T_i $
    
    with finite support, where each $T_i$ denotes the forcing we obtain if we force with the nodes of the tree  as conditions using the tree order as the partial order, then in the resulting generic extension $V[G]$, for every $\alpha \notin A$, $V[G] \models `` T_{\alpha}$ is a Suslin tree$"$. 
    
    \begin{theorem}\label{DefinitionIndependentSequence}(see \cite{Separation without reduction})
    
    Assume that $M$ is an uncountable, transitive model of $\ZF^{-}+``\aleph_1$ exists$"$ satisfying $\aleph_1^M=\aleph_1^L$. Then there is an independent sequence $\vec{S}=(S_\alpha\mid\alpha<\omega_1)$ of Suslin trees in $L$, and the sequence $\vec{S}$ is uniformly $\Sigma_1(\{\omega_1\})$-definable over $M$. To be more precise, there is a $\Sigma_1$-formula $\phi$  with $\omega_1$ as the unique parameter, which does not depend on the model $M$, such that the relation $\{ (t,\gamma) \mid \gamma < \omega_1 \land t \in T^{\gamma}  \}$, where $T^{\gamma}$ denotes the $\gamma$-th level of $T$, is definable over $M$ using $\phi$. 
    
    \end{theorem}
The trees from $\vec{S}$ will later be used to define several related coding forcings that are the main tool in our proofs.

We briefly introduce the almost disjoint coding forcing due to R. Jensen and R. Solovay (\cite{JensenSolovay}). We will identify subsets of $\omega$ with their characteristic function and will use the word reals for elements of $2^{\omega}$ and subsets of $\omega$ respectively.
    
    Let $D=\{d_{\alpha} \, : \, \alpha < \aleph_1 \}$ be a family of almost disjoint subsets of $\omega$, i.e. a family such that if $r, s \in D$ then 
    
    $r \cap s$ is finite. Let $X\subseteq\omega_1$ be a set of ordinals. Then there 
    
    is a ccc forcing, the almost disjoint coding $\mathbb{A}_D(X)$ which adds 
    
    a new real $x$ which codes $X$ relative to the family $D$ in the following way
    
    $$\alpha \in X \text{ if and only if } x \cap d_{\alpha} \text{ is finite.}$$
    
    \begin{definition}\label{definitionadcoding}
    
     The almost disjoint coding $\mathbb{A}_D(X)$ relative to an almost disjoint family $D$ consists of
    
     conditions $(r, R) \in [\omega]^{<\omega} \times D^{<\omega}$ and
    
     $(s,S) < (r,R)$ holds if and only if
    
     \begin{enumerate}
    
      \item $r \subset s$ and $R \subset S$.
    
      \item If $\alpha \in X$ and $d_{\alpha} \in R$ then $r \cap d_{\alpha} = s \cap d_{\alpha}$.
    
     \end{enumerate}
    
    \end{definition}
We shall briefly discuss the $L$-definable, $\aleph_1^L$-sized almost disjoint family of reals $D$ we will use throughout this article. The family $D$ is the canonical almost disjoint family one obtains by recursively choosing the $<_L$-least real $x_\beta$ not yet chosen and replacing it with $d_\beta\subseteq\omega$, where $d_\beta$ codes the initial segments of $x_\beta$ using fixed recursive bijections between $\omega$ and $\omega^{<\omega}$. The definition of $D$ is uniform over every uncountable transitive $\ZFP$-model $M$ satisfying $\aleph_1^M=\aleph_1^L$, since inside $M$ we can compute $L$ correctly up to $\aleph_1^L$ and then apply the above definition in $M$'s version of $L$. Even more is true: if $M$ is a countable transitive model of $\ZFP+``\aleph_1$ exists and $\aleph_1=\aleph_1^L$'', then $M$ computes $D\upharpoonright\omega_1^M$ correctly. Again, $M$ can correctly define an initial segment of $L$ which suffices to calculate $D\upharpoonright\omega_1^M$.

    \section{Coding machinery}
    
    We continue with the construction of the appropriate notions of forcing which we want to use in our proof. The goal is to first define a coding forcings $\operatorname{Code} (x)$ for reals $x$, which will force for $x$ that a certain $\Sigma^1_3$-formula $\Phi(x)$ becomes true in the resulting generic extension. The coding method is basically the same as in \cite{Reduction} and \cite{Uniformization} and literally the same as in \cite{Separation without reduction}.
   
    In a first step we add $\omega_1$-many $\omega_1$-Cohen subsets with a countably supported product, $$\forceP^1:= \prod_{\alpha< \omega_1} \mathbb{C} (\omega_1).$$
    
    Note that this forcing is itself $\sigma$-closed so no reals are added and $\vec{S}$ is still an independent sequence of Suslin trees. In a second step we force over $L^{\forceP^1}$ to destroy all members of $\vec{S}$ via generically adding an $\omega_1$-branch, that is we  form $\forceP^0:=\prod_{\alpha \in \omega_1} S_{\alpha}$ with finite support. Note that this is an $\aleph_1$-sized, ccc forcing over $L$ and also $L^{\forceP^1}$, and the forcing is independent of the actual model in which it is computed. The two step iteration can be thus conceived as a product of two factors $(\prod_{i < \omega_1} \mathbb{C} (\omega_1))^L \times \prod_{\alpha \in \omega_1} S_{\alpha}$.  In the generic extension $\aleph_1$ is preserved and $\CH$ remains to be true.  
   
     We use $W$ to denote this generic extension of $L$, that is we let $g^0 \times g^1$ be a $\forceP^0 \times \forceP^1$ generic filter over $L$  where $\forceP^0$ adds cofinal branches to each member of $\vec{S}$ and $\forceP^1$ adds $\aleph_1$-many $\omega_1$-Cohen subsets, then \[W=L[g^0 \times g^1].\]

     Let $x \in W$ be a real, let $m,k \in \omega$, and let $\eta<\omega_1$. We write $(m,k,1,x)$ for a real coding this quadruple in a recursive way. The forcing $\operatorname{Code}(m,k,1,x,\eta)$ codes $(m,k,1,x)$ into $\vec{S}$ by almost disjoint coding a set $Y\subseteq\omega_1$:
     \[
     \operatorname{Code}(m,k,1,x,\eta):=\mathbb A_D(Y).
     \]
     The forcing $\operatorname{Code}(m,k,0,x,\eta)$ is defined in the same way. The parameter $\eta$ selects the coordinate $g^1_\eta$ used to define $h$; it is not part of the quadruple written into $\vec S$. We now define $Y$.

     To ease notation we let $g\subseteq\omega_1$ be $g^1_\eta$, where $g^1_\eta$ is the $\eta$-th coordinate of the $\prod_{\alpha<\omega_1}\mathbb C(\omega_1)$-generic filter over $L[g^0]$. We let $\rho:([\omega_1]^\omega)^L\to\omega_1$ be some canonically definable, constructible bijection between these two sets. We use $\rho$ and $g$ to define the set $h\subseteq\omega_1$, which will be the set of indices of $\omega$-blocks of $\vec S$ where we code the characteristic function of $(m,k,1,x)$. Let
     \[
     h:=\{\rho(g\cap\alpha):\alpha<\omega_1\}
     \]
     and let
     \begin{align*}
     A:=&\{\omega\gamma+2n\mid\gamma\in h,\ n\notin(m,k,1,x)\}\\
       &\cup\{\omega\gamma+2n+1\mid\gamma\in h,\ n\in(m,k,1,x)\}.
     \end{align*}

     Let
     \[
     G_A:=g^0\upharpoonright A.
     \]
     Then $G_A$ is the generic filter for the finite-support product $\prod_{\beta\in A}S_\beta$. Choose $X\subseteq\omega_1$ so that $X$ codes $G_A$ and
     \[
     L[X]=L[G_A].
     \]
     Thus $X$ codes no further information. The set $A$ is recovered from the coordinates of $G_A$.

     If $\gamma\in h$, the branch pattern in the $\omega$-block of $\vec S$ starting at $\gamma$ determines $(m,k,1,x)$:
     \begin{itemize}
     \item[$(\ast)_\gamma$] $n\in(m,k,1,x)$ if and only if $S_{\omega\gamma+2n+1}$ has an $\omega_1$-branch, and $n\notin(m,k,1,x)$ if and only if $S_{\omega\gamma+2n}$ has an $\omega_1$-branch.
     \end{itemize}
     The filter $G_A$ adds a branch through each $S_\beta$ with $\beta\in A$. If $\xi\notin A$, then $S_\xi$ remains Suslin in $L[G_A]=L[X]$ by the independence of $\vec S$. Hence $(\ast)_\gamma$ holds.

     We note that we can apply an argument resembling David's trick (\cite{David}) in this situation. We rewrite the information of $X\subseteq\omega_1$ as a subset $Y\subseteq\omega_1$ using the following line of reasoning. Keeping Theorem~\ref{DefinitionIndependentSequence} in mind, any transitive, $\aleph_1$-sized model of $\ZFP$ containing $X$ defines $\vec S$ correctly and decodes from $X$ the pattern $(\ast)_\gamma$ for every $\gamma\in h$. Choose a sufficiently large $\theta$ such that $M=L_\theta[X]$ is such a model, and let $X_M\subseteq\omega_1$ code the pointed structure $(M,\in,X)$. Then, for any uncountable $\beta$ such that $L_\beta[X_M]\models\ZFP$,
     \[
     L_\beta[X_M]\models\text{``The model decoded out of $X_M$ satisfies $(\ast)_\gamma$ for every $\gamma\in h$.''}
     \]
     In particular there will be an $\aleph_1$-sized ordinal $\beta$ as above and we can fix a club $C\subseteq\omega_1$ and a sequence $(M_\alpha:\alpha\in C)$ of countable elementary submodels of $L_\beta[X_M]$ such that
     \[
     \forall\alpha\in C\,(M_\alpha\prec L_\beta[X_M]\land M_\alpha\cap\omega_1=\alpha).
     \]
     Now let the set $Y\subseteq\omega_1$ code the pair $(C,X_M)$ such that the odd entries of $Y$ code $X_M$ and, if $E(Y)$ denotes the set of even entries of $Y$ and $\{c_\alpha:\alpha<\omega_1\}$ is the enumeration of $C$, then
     \begin{enumerate}
     \item $E(Y)\cap\omega$ codes a well-ordering of type $c_0$.
     \item $E(Y)\cap[\omega,c_0)=\emptyset$.
     \item For all $\beta$, $E(Y)\cap[c_\beta,c_\beta+\omega)$ codes a well-ordering of type $c_{\beta+1}$.
     \item For all $\beta$, $E(Y)\cap[c_\beta+\omega,c_{\beta+1})=\emptyset$.
     \end{enumerate}

     We make these choices inside $L[X]$. Thus $Y\in L[X]$, while $Y$ decodes $X_M$ and hence $X$. Consequently,
     \[
     L[Y]=L[X]=L[G_A].
     \]

     We obtain
     \begin{itemize}
     \item[$({\ast}{\ast})$] For any countable transitive model $M$ of ``$\ZFP$ and $\aleph_1$ exists'' such that $\omega_1^M=(\omega_1^L)^M$ and $Y\cap\omega_1^M\in M$, the model $L[Y\cap\omega_1^M]^M$ contains an $\aleph_1^M$-sized transitive model $N$ which satisfies $(\ast)_\gamma$ for $(m,k,1,x)$ and every $\gamma\in h\cap\omega_1^M$.
     \end{itemize}

     We use the existential form in $({\ast}{\ast})$. After $Y$ is coded by a real $r$, the conclusion holds for every transitive model with the same $\omega_1$ that contains $r$, as noted below.

     Thus we have a local version of the property $(\ast)$. We have finally defined the desired set $Y$ and now use
     \[
     \operatorname{Code}(m,k,1,x,\eta):=\mathbb A_D(Y)
     \]
     relative to our previously defined almost disjoint family $D\in L$ (see the paragraph after Definition~\ref{definitionadcoding}) to code the set $Y$ into one real $r$. This forcing only depends on the subset of $\omega_1$ we code, so $\mathbb A_D(Y)$ is independent of the surrounding universe in which we define it, as long as it has the right $\omega_1$ and contains $Y$. The effect of the coding forcing $\operatorname{Code}(m,k,1,x,\eta)$ is that it generically adds a real $r$ such that
     \begin{itemize}
     \item[$({\ast}{\ast}{\ast})$] For any countable transitive model $M$ of ``$\ZFP$ and $\aleph_1$ exists'', such that $\omega_1^M=(\omega_1^L)^M$ and $r\in M$, $M$ can construct its version of $L[r]$, denoted by $L[r]^M$, which in turn thinks that there is a transitive $\ZFP$-model $N$ of size $\aleph_1^M$ such that $N$ believes $(\ast)$ for $(m,k,1,x)$ and every $\gamma\in h\cap\omega_1^M$.
     \end{itemize}

     Indeed, if $r$ and $M$ are as above, then $M$ and $L[r]^M$ compute $D\upharpoonright\omega_1^M$ correctly, as discussed below Definition~\ref{definitionadcoding}. As a consequence, $L[r]^M$ contains the set $Y\cap\omega_1^M$, where $Y\subseteq\omega_1$ is as in $({\ast}{\ast})$. So in $L[Y\cap\omega_1^M]^M$, there is an $\aleph_1^M$-sized transitive $N$ which models $(\ast)_\gamma$ for every $\gamma\in h\cap\omega_1^M$, as claimed.

     Note that $({\ast}{\ast}{\ast})$ is a $\Pi^1_2$-formula in the parameters $r$ and $(m,k,1,x)$, since $r$ codes the set $h\cap\omega_1^M$. We say that $(m,k,1,x)$ is \emph{written into $\vec S$}, or \emph{coded into $\vec S$}, if there is a real $r$ such that $({\ast}{\ast}{\ast})$ holds for $r$ and $(m,k,1,x)$. Thus, when applied over $W$, the forcing $\operatorname{Code}(m,k,1,x,\eta)$ adds such a real $r$.

     The forcing $\operatorname{Code}(m,k,0,x,\eta)$ codes $(m,k,0,x)$ into $\vec S$ in the same way. The full branch pattern remains correct in $L[r]$. Indeed, $L[Y]=L[X]=L[G_A]$, so in $L[Y]$ every tree indexed by $A$ has a branch and every tree not indexed by $A$ is Suslin. The forcing $\mathbb A_D(Y)$ is $\sigma$-centered and therefore preserves Suslin trees. Since $r$ is $\mathbb A_D(Y)$-generic over $L[Y]$ and $Y\in L[r]$, we have $L[r]=L[Y][r]$. Hence, for every $\gamma\in h$ and $n\in\omega$,
     \begin{align*}
     n\in(m,k,1,x)&\Longrightarrow L[r]\models ``S_{\omega\gamma+2n+1}\text{ has an $\omega_1$-branch and }S_{\omega\gamma+2n}\text{ is Suslin}'',\\
     n\notin(m,k,1,x)&\Longrightarrow L[r]\models ``S_{\omega\gamma+2n}\text{ has an $\omega_1$-branch and }S_{\omega\gamma+2n+1}\text{ is Suslin}''.
     \end{align*}

    \section{Allowable forcings}

    Next we define the set of forcings which we will use in our proof. We aim to iterate the coding forcings we defined in the last section. We need to generalize the notion of allowable, introduced in \cite{Reduction}, as  we  want to force towards 
\begin{itemize}
\item $\Pi^1_3$-reduction,
\item the failure of $\Pi^1_3$-uniformization,
\item $\Sigma^1_4$-uniformization and
\item a good $\Sigma^1_5$-well-order.
\end{itemize}
These different tasks will be tackled using iterated coding forcings which obey the following definition.
    
    \begin{definition}\label{0-allowable forcing}

    Let $W$ be our ground model and let $\alpha<\omega_1$. We take the bookkeeping function $F$ to be an element of $L$ with domain $\alpha$. Since $W$ is a set-generic extension of $L$, the values of $F$ are understood as $L$-names for the tuples of names and objects used below; we suppress this outer layer of names throughout.

    A finite support iteration $\forceP=((\forceP_\beta,\dot{\forceQ}_\beta):\beta<\alpha)$ is called allowable (relative to the bookkeeping function $F$) if $F$ determines $\forceP$ inductively as follows. Suppose that $\beta<\alpha$, that $\forceP_\beta$ has been defined, and let $G_\beta$ be a $\forceP_\beta$-generic filter over $W$. For $\gamma<\beta$, write $G_\gamma=G_\beta\upharpoonright\forceP_\gamma$. We say that $\eta<\omega_1$ is \emph{free at stage $\beta$} if there are no $\gamma<\beta$ and no finite tuple $\vec z$ such that
    \[
    \dot{\forceQ}_\gamma^{G_\gamma}=\operatorname{Code}(\vec z,\eta).
    \]
    Otherwise, $\eta$ has already been used for coding. We also fix an arity-sensitive recursive coding of finite tuples. In particular, the tagged code of the singleton tuple $(x)$ is distinct from the tagged code of the pair $(x,y)$. We define the next factor according to the form of $F(\beta)$.

    \begin{itemize}
    \item[(1)] Assume that $F(\beta)=(\dot m,\dot k,\dot l,\dot x,\dot\eta)$, where these are $\forceP_\beta$-names. Let $\dot x^{G_\beta}=x$, $\dot m^{G_\beta}=m$, $\dot k^{G_\beta}=k$, $\dot l^{G_\beta}=l$, and $\dot\eta^{G_\beta}=\eta$. We assume that $x$ is a real, that $m$ and $k$ are natural numbers coding two $\Pi^1_3$-sets $A_m$ and $A_k$, that $l\in\{0,1\}$, and that $\eta<\omega_1$.

    If $\eta$ is not free, we let $\dot{\forceQ}_\beta^{G_\beta}$ be the trivial forcing. If $\eta$ is free, we let
    \[
    \dot{\forceQ}_\beta^{G_\beta}=\forceP(\beta)^{G_\beta}:=\operatorname{Code}(m,k,l,x,\eta).
    \]

    \item[(2)] If, on the other hand, $F(\beta)=(\dot m,\dot x,\dot y,\dot\eta)$ and $\dot m$ is the G\"odel number of a $\Pi^1_3$-set in the plane, $\dot x,\dot y$ are both $\forceP_\beta$-names of reals, and $\dot\eta$ is a name for a countable ordinal, then we define
    \[
    \dot{\forceQ}_\beta^{G_\beta}:=\operatorname{Code}(0,0,\dot x^{G_\beta},\dot y^{G_\beta},\dot\eta^{G_\beta})
    \]
    provided $\eta$ is free. If $\eta$ is not free, we let $\eta'$ be the least ordinal which is free and use
    \[
    \dot{\forceQ}_\beta^{G_\beta}:=\operatorname{Code}(0,0,\dot x^{G_\beta},\dot y^{G_\beta},\eta')
    \]
    instead. To avoid an ambiguity in the definition we also declare that $0$ is not the G\"odel number of any formula. The seemingly redundant information when coding $(0,0,x,y)$ instead of just $(0,x,y)$ avoids an ambiguous definition again.

    \item[(3)] If $F(\beta)=(\dot m,\dot x,\dot y,\dot a_0,\dot a_1,\dot\eta)$, where $\dot m^{G_\beta}=m\in\omega$ and $m$ is the G\"odel number of a $\Sigma^1_4$-formula, then we use $\operatorname{Code}(m,x,y,a_0,a_1,\eta)$, provided $\eta$ is free. Otherwise we use the least $\eta'$ which is free and code there.

    \item[(4)] If $F(\beta)=(\dot m,\dot x,\dot y,\dot\eta)$ and $\dot x,\dot y$ are both names of reals whereas $\dot m$ is the G\"odel number of a $\Sigma^1_5$-formula, then we use $\operatorname{Code}(1,1,x,y,\eta)$, provided $\eta$ is free, and use the least $\eta'$ which is free otherwise for coding $(1,1,x,y)$ there. Again to avoid ambiguity we assume that $1$ is not the G\"odel number of a formula and code $(1,1,x,y)$ instead of just $(1,x,y)$.

    \item[(5)] Finally, if $F(\beta)$ is of the form $(\dot n,\dot x,\dot\eta)$, for $\dot x$ a $\forceP_\beta$-name of a real, $\dot n$ a $\forceP_\beta$-name of a natural number and $\dot\eta$ the $\forceP_\beta$-name of an ordinal below $\omega_1$, then we let
    \[
    \dot{\forceQ}_\beta^{G_\beta}=\forceP(\beta):=\operatorname{Code}(n,x,\eta)
    \]
    provided $\eta$ is free for coding. Otherwise we let $\eta'$ be the least ordinal which is free and code at $\eta'$ instead. This item also explains why we coded $(0,0,x,y)$ in item~(2) and $(1,1,x,y)$ in item~(4).
    \end{itemize}
    \end{definition}
 We add as a mildly clarifying remark that the forcings from (1) (2) and (5) are used to work towards a universe where $\Pi^1_3$-reduction holds and $\Pi^1_3$-uniformization fails; and the forcings from (3) and (4) are used to force $\Sigma^1_4$-uniformization and a good $\Sigma^1_5$-well-order of the reals respectively.

Allowable forcings form the base class of the thinning-out process. We also call them $0$-allowable and write $\mathcal A_0$ for their class.

    Every allowable forcing $\forceP$ can be written over $L$ as\footnote{Here and later we will simply write $\Asterisk \forceP(\alpha)$ for the forcing iteration one obtains when using the $\forceP(\alpha)$'s as factors. More precisely, if $(\forceP_{\alpha},\dot{\forceQ}_{\alpha}\mid\alpha<\kappa)$ is the usual notation for an iteration of length $\kappa$, then $\Asterisk_{\alpha<\kappa}\dot{\forceQ}_{\alpha}$ denotes $\forceP_{\kappa}$.}
\[
\forceP^1\ast\forceP^0\ast
\Bigl(\Asterisk_{\xi<\delta}\forceP(\xi)\Bigr),
\]
where, consistently with Section~3,
\[
\forceP^1=\prod_{\xi<\omega_1}\mathbb C(\omega_1),
\qquad
\forceP^0=\prod_{\xi<\omega_1}S_\xi,
\]
and the factors $\forceP(\xi)$ are almost disjoint coding forcings $\mathbb A_D(Y_\xi)$.

We note that $\forceP^0$ and $\Asterisk_{\xi<\delta}\forceP(\xi)$ are ccc. Every almost disjoint coding factor in the allowable iteration uses exactly one coordinate of the $\forceP^1$-generic as its coding area. Since the allowable iteration has countable length, there is a countable ordinal $\alpha<\omega_1$ such that all coding areas used by $\forceP$ are contained in $\alpha+1$. Thus we may factor
\[
\forceP^1=
\prod_{i\leq\alpha}\mathbb C(\omega_1)
\times
\prod_{i>\alpha}\mathbb C(\omega_1)
\]
and rearrange $\forceP$ as
\[
\forceP=
\Bigl(\prod_{i\leq\alpha}\mathbb C(\omega_1)\Bigr)
\ast
\left(
  \left(\forceP^0\ast
  \Asterisk_{\xi<\delta}\forceP(\xi)\right)
  \times
  \prod_{i>\alpha}\mathbb C(\omega_1)
\right).
\]

Work now over
\[
V=L^{\prod_{i\leq\alpha}\mathbb C(\omega_1)}.
\]
The forcing $\prod_{i>\alpha}\mathbb C(\omega_1)$ is $\sigma$-closed in $V$, while
\[
\forceP^0\ast\Asterisk_{\xi<\delta}\forceP(\xi)
\]
is ccc. By Easton's Lemma (see Lemma~15.19 of \cite{Jech}), the Cohen tail remains $\omega$-distributive after forcing with this ccc factor. Consequently every real in $L^{\forceP}$ already belongs to
\[
L^{\left(\prod_{i\leq\alpha}\mathbb C(\omega_1)\right)
  \ast\left(\forceP^0\ast\Asterisk_{\xi<\delta}\forceP(\xi)\right)}.
\]
In particular, every name for a real arising from an allowable forcing can be taken to depend on only countably many coding areas.

As a second and related remark, every allowable forcing $\forceP\in W$ is already correctly definable over a proper inner model of $W$. The bookkeeping function $F$ belongs to $L$. Let $I_0,J_0\subseteq\omega_1$ be countable sets of coordinates which suffice to evaluate the names occurring in $F\upharpoonright\alpha$. A nontrivial factor $\operatorname{Code}(\vec x,\eta)$ also uses the branches $g^0\upharpoonright A_\eta$, where $A_\eta\subseteq\omega_1$ has size $\aleph_1$.

Let $K\subseteq\omega_1$ be the countable set of coding areas used by $\forceP$, and, for $\eta\in K$, let $A_\eta$ denote the set $A$ used in the factor with coding area $\eta$. Set
\[
J:=J_0\cup K
\qquad\text{and}\qquad
I:=I_0\cup\bigcup_{\eta\in K}A_\eta.
\]
Then $J$ is countable and $\forceP$ is correctly defined in
\[
L[g^0\upharpoonright I,g^1\upharpoonright J].
\]
If $K\ne\emptyset$, then $I$ has size $\aleph_1$. It is also co-uncountable. Indeed, choose a coding area $\xi\notin K$. For every $\eta\in K$, the sets $h_\xi$ and $h_\eta$ have countable intersection. Hence $h_\xi\setminus\bigcup_{\eta\in K}h_\eta$ has size $\aleph_1$. The corresponding $\omega$-blocks are disjoint from $\bigcup_{\eta\in K}A_\eta$, and only countably many of their tree coordinates belong to $I_0$. Thus $\omega_1\setminus I$ has size $\aleph_1$. If $K=\emptyset$, enlarge $I_0$ to an uncountable, co-uncountable set and call the result $I$. In either case, both $I$ and $\omega_1\setminus I$ have size $\aleph_1$, while $J$ is countable, and
\[
L[g^0\upharpoonright I,g^1\upharpoonright J]\subsetneq W.
\]
In particular, there are $\aleph_1$-many trees from $\vec S$ and $\aleph_1$-many coding areas which do not enter the definition of $\forceP$.
        If $\forceP \in W$ is a forcing such that there is an $\alpha < \omega_1$ and an $F \in W$, $F: \alpha \rightarrow W$ such that $\forceP$ is allowable with respect to $F$, then we often just drop the $F$ and simply say that $\forceP \in W$ is allowable.
    As mentioned already informally, every allowable forcing uniquely defines a countable set of coding areas it uses with its coding forcings. 
    
    \begin{definition}
    
    Let $\forceP=((\forceP_{\alpha},\dot{\forceQ}_{\alpha}=\forceP({\alpha}))\mid\alpha<\delta)$ be an allowable forcing, and let $G\subset\forceP$ be generic over $W$. We define
    \[
    C^G:=\left\{\eta<\omega_1:\exists\beta<\delta\,\exists\vec z\,
    \bigl(\dot{\forceQ}_{\beta}^{G_{\beta}}=\operatorname{Code}(\vec z,\eta)\bigr)\right\},
    \]
    where $\vec z$ ranges over finite tuples of the appropriate arity for one of the coding forcings from Definition~\ref{0-allowable forcing}. Thus $C^G$ records exactly the coding areas which occur in actual nontrivial coding factors of the iteration, independently of which clause of Definition~\ref{0-allowable forcing} produced the factor.
    
    We also let
    \[
    C^{\forceP}:=\{\eta<\omega_1:\exists p\in\forceP\ (p\Vdash\eta\in C^{\dot G})\}.
    \]
    
    \end{definition}
    
    As noted already above, $C^G$ and $C^{\forceP}$ are countable for every allowable $\forceP$. In particular, all later disjointness requirements of the form $C^{\forceP}\cap C^{\forceQ}=\emptyset$ refer to the complete sets of coding areas used by the two forcings, including coding factors arising from any of clauses~(1)--(5) of Definition~\ref{0-allowable forcing}. We derive some very easy properties of allowable forcings.
    
     \begin{lemma}(see \cite{Separation without reduction})

    \begin{enumerate}
    
    \item If $\forceP=((\forceP(\beta), \dot{\forceQ}_{\beta}) \, : \, \beta < \delta) \in W$ is allowable, then for every $\beta<\delta$,
    \[
    \forceP_\beta\Vdash |\dot{\forceQ}_\beta|\leq\aleph_1.
    \]
    Every nontrivial coding factor is forced to have size $\aleph_1$.
    
    \item Every allowable forcing over $W$ is ccc and thus preserves cardinals.
    
    \item Every allowable forcing over $W$ preserves $\CH$. Furthermore, if $\forceP= (\forceP_{\alpha},\dot{\forceQ}_{\alpha}) \, : \, \alpha < \omega_1) \in W$ is an $\omega_1$-length iteration such that each initial segment of the iteration is allowable over $W$, then $W^{\forceP} \models \CH$.
    
    \item The product of two allowable forcings $\forceP$ and $\forceQ$ can be densely embedded into an allowable forcing provided that $C^{\forceP} \cap C^{\forceQ}=\emptyset$.
    
    \end{enumerate}
    
    \end{lemma}

    Let $\forceP=((\forceP_\beta,\dot{\forceQ}_\beta)\mid\beta<\delta)$ be allowable and let $G\subseteq\forceP$ be generic over $W$. We say that a finite tuple $\vec z$ is \emph{coded by $\forceP$} if for some $\beta<\delta$ and some $\eta<\omega_1$,
    \[
    \dot{\forceQ}_\beta^{G_\beta}=\operatorname{Code}(\vec z,\eta).
    \]
    Thus the objects coded by an allowable forcing are exactly the tuples occurring in its actual nontrivial coding factors; they need not coincide with all reals or tuples mentioned by the bookkeeping function.

    For a finite tuple $\vec z$, let $\Phi(\vec z)$ denote the corresponding $\Sigma^1_3$ coding statement from Section~3, namely that there is a real $r$ witnessing that $\vec z$ is coded into $\vec S$.

    A crucial property of allowable forcings is that no additional tuples are coded accidentally.

    \begin{lemma}\label{nounwantedcodes}(see \cite{Reduction}, \cite{Uniformization})
    Let $\forceP=((\forceP_\beta,\dot{\forceQ}_\beta)\mid\beta<\delta)\in W$ be allowable and let $G\subseteq\forceP$ be generic over $W$. Then, in $W[G]$, for every finite tuple $\vec z$,
    \[
    \Phi(\vec z)
    \quad\Longleftrightarrow\quad
    \exists\beta<\delta\,\exists\eta<\omega_1\,
    \bigl(\dot{\forceQ}_\beta^{G_\beta}=\operatorname{Code}(\vec z,\eta)\bigr).
    \]
    \end{lemma}

    \section{Thinning out}

We define next a derivative of the class of allowable forcings. These derivatives can be applied transfinitely often, yielding smaller and smaller nonempty classes of allowable forcings. We write $\mathcal A_0$ for the class of $0$-allowable forcings. Eventually the derivative operator acts on a nonempty class which cannot be thinned out further; this class is a fixed point of the operation. We denote the fixed point by $\mathcal A_\infty$, and its members are called $\infty$-allowable. They are the coding forcings which we use to force $\Pi^1_3$-reduction. We emphasize that the other tasks, namely forcing $\Sigma^1_n$-uniformization for $n\geq 4$ and the failure of $\Pi^1_3$-uniformization, play no role in the definition of this thinning-out process. These tasks will be incorporated later, after the thinning-out operator has been defined.

As there will be several longer definitions, we first motivate the thinning-out process and provide the intuition behind the construction. A more detailed discussion can be found in \cite{Reduction}.

\subsection{Informal discussion of the idea}

We proceed with an informal discussion of the main ideas of the proof. We consider an arbitrary pair $A_m$ and $A_k$ of $\Pi^1_3$-sets and want to find reducing sets $D^0_{m,k}$ and $D^1_{m,k}$, that is, sets with the following properties:
\begin{enumerate}
\item $D^0_{m,k}$ and $D^1_{m,k}$ are both $\Pi^1_3$;
\item $D^0_{m,k}\subseteq A_m$ and $D^1_{m,k}\subseteq A_k$;
\item $D^0_{m,k}\cap D^1_{m,k}=\emptyset$;
\item $D^0_{m,k}\cup D^1_{m,k}=A_m\cup A_k$.
\end{enumerate}

The ansatz is to use the two types of coding forcings $\operatorname{Code}(m,k,0,x,\eta)$ and $\operatorname{Code}(m,k,1,x,\eta)$ to define reducing sets for the fixed pair of $\Pi^1_3$-formulas $\varphi_m$ and $\varphi_k$. The two candidates for reducing sets are defined by the $\Pi^1_3$-formulas
\begin{align*}
D^0_{m,k}:={}&\{x:(m,k,0,x)\text{ is not coded into }\vec S\},
\end{align*}
and
\begin{align*}
D^1_{m,k}:={}&\{x:(m,k,1,x)\text{ is not coded into }\vec S\}.
\end{align*}
Here, as in Section~3, saying that $(m,k,i,x)$ is coded into $\vec S$ means that there is a real $r$ such that $({\ast}{\ast}{\ast})$ holds for $r$ and $(m,k,i,x)$.

As always there will be a bookkeeping function $F$ which hands us at every stage $\beta<\omega_1$ a name for a real $x$, and the task is to decide which one of the two forcings $\operatorname{Code}(m,k,0,x,\eta)$ and $\operatorname{Code}(m,k,1,x,\eta)$ we want to use at that stage. In other words, we decide at stage $\beta$ whether we place $x$ into $D^0_{m,k}$ or into $D^1_{m,k}$.

We approach this problem as follows. As a first observation, if a real $x$ cannot be forced out of $A_m$ with a forcing in $\mathcal A_0$, then we can safely put it into $D^0_{m,k}$ by using $\operatorname{Code}(m,k,1,x,\eta)$ for some $\eta$ and ensure that we never put $x$ into $D^1_{m,k}$ later. If $x$ can be forced out of $A_m$ with a forcing in $\mathcal A_0$ but cannot be forced out of $A_k$ with such a forcing, then it is safe to place $x$ into $D^1_{m,k}$. In the remaining case, $x$ can be forced out of both $A_m$ and $A_k$, and we let the bookkeeping function decide where to place $x$.

The class obtained from this rule is a first approximation to a class of iterations which should solve the problem of finding reducing sets $D^0_{m,k}$ and $D^1_{m,k}$. We denote it by $\mathcal A_1$ and call its members $1$-allowable. We note that when forcing with a member of $\mathcal A_1$, we ask the wrong questions when trying to place $x$. Indeed, if we run a $1$-allowable iteration, then whenever we ask whether some $x$ can be forced out of $A_m$ with a forcing in $\mathcal A_0$, we should instead ask whether $x$ can be forced out of $A_m$ with a forcing in $\mathcal A_1$, since this is the class to which our iteration belongs.

Thus it is reasonable to add a second question at each stage of a $1$-allowable forcing. Given a real $x$, we first ask whether $x$ can be forced out of $A_m$ or $A_k$ with a forcing in $\mathcal A_0$. If the answer is yes, we then ask whether $x$ can be forced out of $A_m$ or $A_k$ with a forcing in $\mathcal A_1$. If the answer to the second question is no, we can safely place $x$ into $D^0_{m,k}$ or $D^1_{m,k}$, since the forcing we are defining will itself be $1$-allowable. The resulting class is denoted by $\mathcal A_2$, and its members are called $2$-allowable. But, looking at the definition of $\mathcal A_2$, we see that we again ask the wrong questions.

These considerations lead to a fixed-point problem. We solve it by repeating the construction transfinitely, obtaining better approximations to the class of iterations needed for the reduction argument. Eventually the sequence stabilizes at the class $\mathcal A_\infty$. Its members are the $\infty$-allowable forcings. They solve the fixed-point problem and will be used to force the reduction property.

\subsection{The derivative operator}

The definition below is uniform over the generic extensions under consideration. Whenever $\mathcal A_\alpha$ is used inside a model $N$, it denotes the class obtained by carrying out the definition in $N$. Thus, at a stage evaluated in $N[G_\beta]$, all quantifiers over $\zeta$-allowable forcings range over the class $\mathcal A_\zeta$ as computed in $N[G_\beta]$.

\begin{definition}\label{AlphaAllowable}
Let $\alpha>0$ and assume that $\mathcal A_\zeta$ has been defined for every $\zeta<\alpha$, uniformly over the generic extensions under consideration. Work in one such model $N$. A finite-support iteration $\forceP$ of length $\delta<\omega_1$ belongs to $\mathcal A_\alpha$ if it is generated by a bookkeeping function $F\in L$ as follows.

Suppose that $\forceP_\beta$ has been defined and let $G_\beta$ be $\forceP_\beta$-generic over $N$. At a reduction stage write
\[
F(\beta)=((\dot x,m,k),(\ell,\dot\eta)),
\]
where $m<k$ code $A_m,A_k$, $\ell\in2$, $x=\dot x^{G_\beta}$, and $\eta=\dot\eta^{G_\beta}<\omega_1$. If $x\in A_m\cup A_k$, choose the least $\zeta<\alpha$ for which, in $N[G_\beta]$, one of the following holds:
\begin{enumerate}
\item[(i)] every $\zeta$-allowable forcing forces $x\in A_m$;
\item[(ii)] every $\zeta$-allowable forcing forces $x\in A_k$.
\end{enumerate}
If both hold, use~(i). In case~(i) set $i=1$, in case~(ii) set $i=0$, and if no such $\zeta$ exists set $i=\ell$. If $x\notin A_m\cup A_k$, also set $i=\ell$. The next factor is
\[
\operatorname{Code}(m,k,i,x,\eta'),
\]
where $\eta'=\eta$ if $\eta$ is free and $\eta'$ is the least free coding area otherwise.

If $F(\beta)$ has one of the forms in clauses~(2)--(5) of Definition~\ref{0-allowable forcing}, use the corresponding rule from that definition. In all other cases use the trivial forcing. At limit stages take the finite-support direct limit.
\end{definition}

A forcing is $\alpha$-allowable if it belongs to $\mathcal A_\alpha$.

\begin{lemma}\label{AllowableDerivativeProperties}
Let $\beta<\alpha$.
\begin{enumerate}
\item $\mathcal A_\alpha$ is definable over $W$, is nonempty, and
\[
\mathcal A_\alpha\subseteq\mathcal A_\beta.
\]
\item If $\forceP,\forceQ\in\mathcal A_\alpha$ and $C^{\forceP}\cap C^{\forceQ}=\emptyset$, then $\forceP\times\forceQ$ densely embeds into a forcing in $\mathcal A_\alpha$.
\end{enumerate}
\end{lemma}

The proof is the induction from \cite{Reduction}. For the product clause, concatenate the bookkeeping functions. The resulting function belongs to $L$.

Modulo isomorphism, there are only set-many partial orders of size at most $\aleph_1$. Hence the decreasing sequence $(\mathcal A_\alpha\mid\alpha\in\operatorname{Ord})$ stabilizes. As the definition of $\alpha$-allowable implies that the set $\mathcal{A}_{\alpha}$ is always nonempty, the decreasing sequence will stabilize at a non-empty set of forcing notions. Let $\theta$ be least such that $\mathcal A_\theta=\mathcal A_{\theta+1}$ and set
\[
\mathcal A_\infty:=\mathcal A_\theta.
\]
Its members are the $\infty$-allowable forcings.
\section{Thinning out while leaving out coding forcings}

One part of the construction is a version of the thinning-out process in which specified coding forcings are omitted. This gives control over the $\Pi^1_3$-relation used later to witness the failure of $\Pi^1_3$-uniformization, while the iteration forces $\Pi^1_3$-reduction.

For a fixed set $B$ of tuples of forcing names, we replace the base class $\mathcal A_0$ by the class of $0$-allowable forcings which avoid the coding forcings associated with $B$, and then apply the derivative from the preceding section. We shall carry out the details now.

Let $B$ be a set of finite tuples of forcing names.

\begin{definition}\label{AvoidingB}
An allowable iteration
\[
\forceP=((\forceP_\beta,\dot{\forceQ}_\beta)\mid\beta<\delta)
\]
\emph{avoids $B$} if, for every $\beta<\delta$ and every $\forceP_\beta$-generic filter $G_\beta$, there are no $\vec x\in B$ and $\eta<\omega_1$ such that
\[
\dot{\forceQ}_\beta^{G_\beta}=\operatorname{Code}(\vec x^{G_\beta},\eta).
\]
\end{definition}

We assume that $B$ contains no tuple $(m,k,i,\dot x)$ used by clause~(1) of Definition~\ref{0-allowable forcing}.

For fixed $B$, set
\[
\mathcal A^B_0:=\{\forceP\in\mathcal A_0\mid\forceP\text{ avoids }B\}.
\]

\begin{definition}\label{AlphaAllowableAvoidingB}
For $\alpha>0$, define $\mathcal A^B_\alpha$ by Definition~\ref{AlphaAllowable}, starting from $\mathcal A^B_0$. Thus the search is over $\zeta<\alpha$, ``$\zeta$-allowable'' means membership in $\mathcal A^B_\zeta$, every factor avoids $B$, and clauses~(2)--(5) of Definition~\ref{0-allowable forcing} are unchanged.
\end{definition}

\begin{lemma}\label{AvoidingBProperties}
Let $\beta<\alpha$.
\begin{enumerate}
\item $\mathcal A^B_\alpha$ is definable, is nonempty, and $\mathcal A^B_\alpha\subseteq\mathcal A^B_\beta$.
\item If $\forceP,\forceQ\in\mathcal A^B_\alpha$ and $C^{\forceP}\cap C^{\forceQ}=\emptyset$, then $\forceP\times\forceQ$ densely embeds into a forcing in $\mathcal A^B_\alpha$.
\end{enumerate}
\end{lemma}

The proof repeats the induction for Lemma~\ref{AllowableDerivativeProperties}, carrying the avoidance condition through each stage. The sequence stabilizes; write
\[
\mathcal A^B_\infty
\]
for its fixed point. The product clause also holds in $\mathcal A^B_\infty$.

\begin{definition}\label{inftallowableAvoidingB}
The members of $\mathcal A^B_\infty$ are the $\infty$-allowable forcings avoiding $B$.
\end{definition}

\subsection{The restart principle}

The fixed-$B$ construction above is the local thinning-out operation needed in the main iteration. It is important that the next operation is only a one-step restart. We do not fix a sequence of forbidden families in advance. Rather, the main construction first uses forcings from its current fixed-point class. If a forcing phase produces new forbidden information, then, in the resulting model, the current class is restricted by the new avoidance requirement and the derivative is run again.

\begin{definition}\label{RestartThinningOut}
Let $M$ be one of the generic extensions under consideration, and compute all classes occurring below in $M$. Suppose that $\mathcal C_\infty$ is the current fixed-point class, as computed in $M$, obtained from a previous application of the derivative from Definition~\ref{AlphaAllowable}. We assume that $\mathcal C_\infty$ is nonempty and has the corresponding product property for forcings which use disjoint coding areas.

Let $D\in M$ be a set of finite tuples of forcing names satisfying the same standing restriction as above: $D$ contains no tuple $(m,k,i,\dot x)$ used by clause~(1) of Definition~\ref{0-allowable forcing}. The set $D$ records only the newly added forbidden information; the restrictions imposed at earlier stages are already incorporated into $\mathcal C_\infty$. Define
\[
\operatorname{Rest}(\mathcal C_\infty,D)_0:=\{\forceP\in\mathcal C_\infty\mid\forceP\text{ avoids }D\}.
\]
For $\alpha>0$, define $\operatorname{Rest}(\mathcal C_\infty,D)_\alpha$ by applying the derivative from Definition~\ref{AlphaAllowable} with $\operatorname{Rest}(\mathcal C_\infty,D)_0$ as the base class. Thus, at derivative rank $\alpha$, the search over $\zeta<\alpha$ is carried out through the classes $\operatorname{Rest}(\mathcal C_\infty,D)_\zeta$, and every factor is required to avoid $D$. Since $\operatorname{Rest}(\mathcal C_\infty,D)_0\subseteq\mathcal C_\infty$, every restriction already incorporated into the current fixed-point class remains in force throughout the restarted thinning-out process.
\end{definition}

\begin{lemma}\label{RestartThinningProperties}
In the setting of Definition~\ref{RestartThinningOut}, let $\beta<\alpha$.
\begin{enumerate}
\item The class $\operatorname{Rest}(\mathcal C_\infty,D)_\alpha$ is uniformly definable over $M$ from $D$ and the defining rule for $\mathcal C_\infty$, is nonempty, and
\[
\operatorname{Rest}(\mathcal C_\infty,D)_\alpha\subseteq\operatorname{Rest}(\mathcal C_\infty,D)_\beta\subseteq\mathcal C_\infty.
\]
In particular, every member of every $\operatorname{Rest}(\mathcal C_\infty,D)_\alpha$ avoids $D$.
\item If $\forceP,\forceQ\in\operatorname{Rest}(\mathcal C_\infty,D)_\alpha$ and
\[
C^{\forceP}\cap C^{\forceQ}=\emptyset,
\]
then $\forceP\times\forceQ$ densely embeds into a forcing in $\operatorname{Rest}(\mathcal C_\infty,D)_\alpha$.
\end{enumerate}
Consequently, the decreasing sequence $(\operatorname{Rest}(\mathcal C_\infty,D)_\alpha\mid\alpha\in\operatorname{Ord})$ stabilizes at a nonempty fixed point, and the product clause holds at the fixed point as well.
\end{lemma}

\begin{proof}
We repeat the induction used for Lemmas~\ref{AllowableDerivativeProperties} and~\ref{AvoidingBProperties}, now with $\operatorname{Rest}(\mathcal C_\infty,D)_0$ as the base class. We indicate why the hypotheses needed for that induction are preserved.

The trivial iteration belongs to $\mathcal C_\infty$ and avoids $D$, so it belongs to $\operatorname{Rest}(\mathcal C_\infty,D)_0$ and, by the same observation at every derivative rank, to each $\operatorname{Rest}(\mathcal C_\infty,D)_\alpha$. Thus all the classes are nonempty. The definition of the derivative is uniform over the generic extensions under consideration, and hence each $\operatorname{Rest}(\mathcal C_\infty,D)_\alpha$ is uniformly definable over $M$ from $D$ and the defining rule for $\mathcal C_\infty$.

The monotonicity argument is the same as in Lemma~\ref{AllowableDerivativeProperties}. Increasing the derivative rank adds further tests against earlier classes and therefore can only remove forcings. Since the restarted process begins with the subclass $\operatorname{Rest}(\mathcal C_\infty,D)_0$ of $\mathcal C_\infty$, we obtain
\[
\operatorname{Rest}(\mathcal C_\infty,D)_\alpha\subseteq\operatorname{Rest}(\mathcal C_\infty,D)_\beta\subseteq\operatorname{Rest}(\mathcal C_\infty,D)_0\subseteq\mathcal C_\infty
\]
whenever $\beta<\alpha$. The requirement that every factor avoid $D$ is carried through the definition at every rank, so every member of these classes avoids $D$.

For the product clause, first consider the base class. If $\forceP,\forceQ\in\operatorname{Rest}(\mathcal C_\infty,D)_0$ use disjoint coding areas, use the product property of $\mathcal C_\infty$. As in the proof of Lemma~\ref{AllowableDerivativeProperties}, the witnessing forcing is obtained by concatenating the two bookkeeping functions. Hence it introduces no coding factor other than a factor already occurring in $\forceP$ or in $\forceQ$. Since both forcings avoid $D$, the resulting forcing avoids $D$ as well and therefore belongs to $\operatorname{Rest}(\mathcal C_\infty,D)_0$.

Now assume inductively that the product clause holds at every rank below $\alpha$, and let $\forceP,\forceQ\in\operatorname{Rest}(\mathcal C_\infty,D)_\alpha$ use disjoint coding areas. Concatenate their bookkeeping functions exactly as in the proof of Lemma~\ref{AllowableDerivativeProperties}. At each stage, the persistence questions appearing in Definition~\ref{AlphaAllowable} are evaluated with respect to the same classes $\operatorname{Rest}(\mathcal C_\infty,D)_\zeta$, $\zeta<\alpha$. The inductive product hypothesis is therefore available whenever the original product argument combines two witness forcings. The concatenated iteration consequently satisfies the rules defining $\operatorname{Rest}(\mathcal C_\infty,D)_\alpha$. Again, no new type of coding factor is introduced, so avoidance of $D$ is preserved. This proves the product clause at rank $\alpha$.

Finally, modulo isomorphism there are only set-many partial orders of size at most $\aleph_1$. Since the sequence is decreasing and every class is nonempty, it stabilizes at a nonempty fixed point. The product clause at the fixed point follows from the clause at the stabilizing rank.
\end{proof}

We write
\[
\operatorname{Rest}(\mathcal C_\infty,D)_\infty
\]
for the fixed point supplied by Lemma~\ref{RestartThinningProperties}, and call the passage from $\mathcal C_\infty$ to $\operatorname{Rest}(\mathcal C_\infty,D)_\infty$ the \emph{restart of the thinning-out process with the new forbidden family $D$}.

No restart is performed when no new forbidden information has been added: in that case the forcing construction simply continues to use the current fixed-point class. In the main construction, a new family $D$ may be defined only after forcing with a member of that class. The forcing stages and the successive restarts must therefore be defined jointly; this stop-and-go recursion will be given where the main iteration is constructed.
\section{Proof of the main theorem}

We will start to define the iteration that will prove the main theorem now. There are four different tasks we have to take care of:
\begin{enumerate}
\item Forcing a failure of $\Pi^1_3$-uniformization,
\item forcing the $\Pi^1_3$-reduction,

\item forcing the $\Sigma^1_4$-uniformization property and
\item forcing  a good $\Sigma^1_5$-wellorder of the reals.
\end{enumerate}

We first produce an intermediate model $W[G_\omega]$ together with a cumulative family $B_\omega$ of forbidden tuples. The role of this preparatory construction is that any later allowable continuation avoiding $B_\omega$ which carries out the required diagonalization at the distinguished reals makes the set $A$ defined below impossible to uniformize by a $\Pi^1_3$-set, as stated in Lemma~\ref{PropertiesOfIntermediateModel}. Over $W[G_\omega]$ we then define the main iteration, interleaving these diagonalization requirements with the other tasks. 

\subsection{Forcing the failure of \texorpdfstring{$\Pi^1_3$}{Pi13}-uniformization}
We shall describe the way how to force a failure of $\Pi^1_3$-uniformization.
The basic first idea is to consider the set
\[
A:= \{ (x,y) \in (\omega^{\omega})^2 \mid (0,0,x,y) \text{ is not coded into $\vec{S}$} \}
\]
which is $\Pi^1_3$.

We will employ our coding forcings to construct a universe where the set $A$ cannot be uniformized by any $\Pi^1_3$-function, via diagonalizing agains all possible $\Pi^1_3$-graphs of functions. This objective is in considerable tension with the simultaneous goal of achieving $\Pi^1_3$-reduction. This tension, however, can be mitigated if we use the coding forcings of the form $\operatorname{Code}(n,x,\eta)$ for $n\in\omega$ from case~(5) in Definition~\ref{0-allowable forcing} of $0$-allowable forcings. These forcings will be used as some sort of yardstick which will mark potentially dangerous places for codes. The idea is inspired by a similar construction in \cite{Separation}.

\subsubsection{The preparatory rules and later continuation}
We first describe the four rules used in the preparatory construction. After the preparatory stage we keep the cumulative family
\[
B_\omega:=\bigcup_{n<\omega}B_n
\]
fixed and explain that under mild restrictions of the further coding forcings a failure of $\Pi^1_3$ uniformization can be forced. These restricitions will be satisfied later in our main iteration so we can argue for a failure of $\Pi^1_3$ uniformization.

First we initialize the forbidden families by setting
\[
B_0:=\emptyset.
\]
During the preparatory step, suppose that we are at a stage $\beta<\omega$, that the iteration up to $\beta$ and a $\forceP_\beta$-generic filter $G_\beta$ have been defined, and that, for some $n\in\omega$, the families $B_j$, $j\leq n$, have been defined. The bookkeeping gives
\[
F(\beta)=(m,x), \text{ where } x\in\omega^\omega \text{ and }m\in\omega.
\]
For an allowable iteration
\[
\forceP=((\forceP_\gamma,\dot{\forceQ}_\gamma)\mid\gamma<\delta),
\]
we say that a coding forcing $\operatorname{Code}(n,\dot R,\dot\eta)$ \emph{occurs as a factor of $\forceP$} if there is a stage $\gamma<\delta$ such that $\dot R$ and $\dot\eta$ are $\forceP_\gamma$-names and
\[
\forceP_\gamma\Vdash
\dot{\forceQ}_\gamma=\operatorname{Code}(n,\dot R,\dot\eta).
\]
We use the same terminology for the other forms of the coding forcing. We now fix the convention for the names which occur in the forbidden families. If $M$ is the current model and $C$ is the cumulative forbidden family already imposed, let
\begin{align*}
    \mathcal N(M,C):=\{(\forceR,\dot z) \mid& \forceR\in M\text{ is a $0$-allowable continuation over $M$} \\&\text{avoiding $C$, and }\dot z\in H(\omega_2)^M\text{ is a nice $\forceR$-name for a real}\}.
\end{align*}

When only the name matters, we suppress the first coordinate and write $\dot z\in\mathcal N(M,C)$. 
We distinguish four cases.

\begin{enumerate}
\item We assume first that $x$ has not been considered before by the bookkeeping function. We also assume that
\begin{align*}
W[G_\beta]\models\exists\forceP\exists y\in\omega^\omega (&\forceP\text{ is $0$-allowable avoiding }\bigcup_{j\leq n}B_j\text{ and }\\
&\text{no forcing of the form }\operatorname{Code}(n,\dot R,\dot\eta)\text{ occurs as a factor of }\forceP\text{ and }\\
&\forceP\Vdash(x,y)\notin A_m).
\end{align*}
In this situation, we first fix the $<$-least such allowable $\forceP$ and use it at stage $\beta$ of our iteration, that is, we let
\[
(\dot{\forceQ}^0_\beta)^{G_\beta}=\forceP.
\]

In a second step we let $R_n$ be a real which codes all the generic reals we have created so far.\footnote{Note that allowable forcings are just a countable length iteration of almost disjoint coding forcings. Hence each allowable extension of $W$ can be written as $W[R]$ for one real $R$ which codes all the a.d. reals added so far.}
Then, if $H^0$ denotes a $\forceP$-generic filter over $W[G_\beta]$, we work over $W[G_\beta][H^0]$ and force with
\[
(\dot{\forceQ}^1_\beta)^{G_\beta}:=\operatorname{Code}(n,R_n,\eta)
\]
for an $\eta$ which is free for coding. Let $M^+$ denote the resulting model after this stage and let
\[
\mathcal N_{n+1}:=\{\dot z\mid \exists\forceR\,((\forceR,\dot z)\in
\mathcal N(M^+,\bigcup_{j\leq n}B_j))\}.
\]
We let
\[
B_{n+1}:=\{(n,\dot z)\mid \dot z\in\mathcal N_{n+1}\}\cup\{(0,0,x,y)\}.
\]
From this point on, we use only $0$-allowable forcings avoiding $\bigcup_{j\leq n+1}B_j$. In other words, we never again use a coding forcing of the form $\operatorname{Code}(n,z,\eta)$ and we also avoid a coding forcing of the form $\operatorname{Code}(0,0,x,y,\eta)$.

\item We assume that case~1 does not apply. In particular, there is no allowable forcing $\forceP$ avoiding $\bigcup_{j\leq n}B_j$ which forces a pair $(x,y)$ out of $A_m$ without containing a factor of the form $\operatorname{Code}(n,\dot R,\dot\eta)$ or $\operatorname{Code}(0,0,x,y,\dot\eta)$. In this situation we decide never again to use a coding forcing of either form. Put
\[
\mathcal N_{n+1}:=\{\dot z\mid \exists\forceR\,((\forceR,\dot z)\in
\mathcal N(W[G_\beta],\bigcup_{j\leq n}B_j))\}.
\]
We set
\[
B_{n+1}:=\{(n,\dot z)\mid\dot z\in\mathcal N_{n+1}\}
\cup
\{(0,0,x,\dot y)\mid\dot y\in\mathcal N_{n+1}\}
\]
and use allowable forcings avoiding $\bigcup_{j\leq n+1}B_j$ from then on. We do not force at this stage.

\item If the real $x$ has already been considered at an earlier stage $\gamma<\beta$, and if case~1 applied at stage $\gamma$ and singled out a real $y$ such that $(0,0,x,y)$ belongs to the cumulative forbidden family, then we pick a real $y'\neq y$ which has not yet been treated and use
\[
\dot{\forceQ}_\beta^{G_\beta}:=\operatorname{Code}(0,0,x,y',\eta)
\]
for some $\eta<\omega_1$ which has not been used for coding.

\item If the real $x$ has been considered earlier and case~2 applied there, then we do not force.
\end{enumerate}

At the end of the preparatory $\omega$-iteration, all families $B_n$, $n<\omega$, have been defined. We then stop this iteration, set
\[
B_\omega:=\bigcup_{n<\omega}B_n,
\]
and call the resulting model $W[G_\omega]$.

For every G\"odel number $m$ of a $\Pi^1_3$-formula in two free variables, let $\beta_m$ be the first preparatory stage at which $m$ is treated, and let $x$ be the real presented at that stage. If case~1 applies at $\beta_m$, set $x_m:=x$, let $y_m$ be the corresponding real singled out by case~1, and call $x_m$ the \emph{distinguished real} associated with $m$. If case~2 applies at $\beta_m$, we assign no distinguished real to $m$.

\begin{lemma}\label{PropertiesOfIntermediateModel}
Let $W[G_\omega]$ and $B_\omega$ be obtained from the preparatory iteration. Let $\widetilde W$ be a suitable outer model of $W[G_\omega]$ obtained by an allowable forcing avoiding $B_\omega$, and let $\forceQ$ be an $\omega_1$-length allowable continuation over $\widetilde W$ which still avoids $B_\omega$. Assume that, for every $m$ for which case~1 applied at the first treatment stage $\beta_m$, the bookkeeping revisits the distinguished real $x_m$ cofinally often and treats every relevant real $y\neq y_m$ according to case~3 above. Then, in every $\forceQ$-generic extension of $\widetilde W$, the set
\[
A:=\{(x,y)\in(\omega^\omega)^2\mid(0,0,x,y)\text{ is not coded into }\vec S\}
\]
is not uniformized by any $\Pi^1_3$-set in the plane.
\end{lemma}

\begin{proof}
Fix a G\"odel number $m$ of a $\Pi^1_3$-formula in two free variables.

Suppose first that case~1 applied at the first treatment stage $\beta_m$. Let $x_m$ be the distinguished real associated with $m$, and let $y_m$ be the real singled out at that stage. The preparatory forcing made $(x_m,y_m)\notin A_m$, and this remains true in all later forcing extensions because the assertion is $\Sigma^1_3$. Moreover, $(0,0,x_m,y_m)$ belongs to $B_\omega$, so no continuation avoiding $B_\omega$ codes this tuple. By the bookkeeping assumption, every real $y'\neq y_m$ which appears in the final extension is eventually treated by case~3 and hence $(0,0,x_m,y')$ is coded. Consequently, in the final extension the $x_m$-section of $A$ is the singleton $\{y_m\}$, while $(x_m,y_m)\notin A_m$. Thus $A_m$ does not uniformize $A$.

Suppose instead that case~2 applied at the first treatment stage $\beta_m$, and let $x$ be the real presented at that stage. Then $B_\omega$ contains $(0,0,x,\dot y)$ for every relevant name $\dot y$ and hence forbids every coding forcing for a tuple $(0,0,x,y)$. The defining failure of case~1 ensures that no allowable continuation avoiding $B_\omega$ can force such a pair $(x,y)$ out of $A_m$. Hence, in the final extension, the $x$-section of both $A$ and $A_m$ contains every real. In particular, $A_m$ is not the graph of a function and therefore does not uniformize $A$.

Thus, for every G\"odel number $m$ of a $\Pi^1_3$-formula in two free variables, $A_m$ does not uniformize $A$. Hence no $\Pi^1_3$-set in the plane uniformizes $A$.
\end{proof}

\subsubsection{The preparatory \texorpdfstring{$\omega$}{omega}-length iteration}\label{PreparatoryOmegaIteration}

The last subsection gave us a mean to force the failure of $\Pi^1_3$-uniformization. As mentioned already, considerable tension arises when trying to combine this with the method to force $\Pi^1_3$-reduction. This is the reason for the addition of the coding forcings $\operatorname{Code} (n,\dot{R},\dot{\eta})$, which are used as the yardstick gauging potentially dangerous coding areas. We shall argue that this yardstick can be made $\Sigma^1_3$ itself, which will unlock a safe way to mitigate this tension.

We start with fixing infinite recursive families of pairwise syntactically distinct $\Pi^1_3$-formulas defining, respectively, the empty set and the whole plane. We use these formulas for padding and arrange that in our preparatory iteration defined above case~1 occurs at every even stage and case~2 at every odd stage. Indeed, if case~1 applies to the next formula to be treated, we treat it at the even stage and use a whole-plane padding formula at the following odd stage. If case~1 fails, we first use an empty-set padding formula. As a consequence in our $\omega$-length preparatory iteration every $\Pi^1_3$-formula is treated while cases~1 and~2 alternate.
At stage $\omega$ we stop, set
\[
B_\omega:=\bigcup_{n<\omega}B_n,
\]
and denote the resulting model by $W[G_\omega]$.

\begin{lemma}\label{DefinabilityR}
For a real $R$, let $(R)_n$ denote the $n$-th part of a recursively definable partition of $R$ into $\omega$-many parts. In $W[G_\omega]$ the following $\Sigma^1_3$-formula holds for exactly one real $R$:
\[
\Psi(R)\equiv
\forall n\in\omega\bigl((2n,(R)_n)\text{ is coded}\bigr)
\land
\forall n\in\omega\setminus\{0\}\bigl((R)_n\notin L[(R)_{n-1}]\bigr).
\]
\end{lemma}
\begin{proof}
We work in $W[G_\omega]$. At every even stage $2n$, case~1 applies. Let $R_{2n}$ be a real coding the countably many reals determining the preparatory iteration up to that point. We then use
\[
\operatorname{Code}(2n,R_{2n},\eta).
\]
For $n=0$, the case~1 test ensures that no factor of the auxiliary forcing used at stage $0$ codes a tuple of the form $(0,z)$, and $B_1$ forbids all such coding forcings afterwards. Hence $R_0$ is the unique real $z$ such that $(0,z)$ is coded.

Now let $n>0$. Every tuple $(2n,z)$ coded before stage $2n$ has $z\in L[R_{2n-2}]$, by the choice of $R_{2n-2}$. The case~1 test ensures that no factor of the auxiliary forcing used at stage $2n$ codes a tuple of the form $(2n,z)$, while $B_{2n+1}$ forbids all such coding forcings afterwards. Moreover,
\[
R_{2n}\notin L[R_{2n-2}],
\]
since $R_{2n}$ also codes the new generic real added by the marker forcing at stage $2n-2$, after $R_{2n-2}$ was chosen. Consequently, $R_{2n}$ is the unique real $z\notin L[R_{2n-2}]$ such that $(2n,z)$ is coded.

Let $R$ be the real determined by
\[
(R)_n=R_{2n}
\qquad(n\in\omega).
\]
Then $\Psi(R)$ holds. Conversely, suppose that $\Psi(R')$ holds. Since $(0,(R')_0)$ is coded, the uniqueness at stage $0$ gives
\[
(R')_0=R_0.
\]
Inductively, suppose that $(R')_{n-1}=R_{2n-2}$ for some $n>0$. Then $\Psi(R')$ implies that $(2n,(R')_n)$ is coded and
\[
(R')_n\notin L[(R')_{n-1}]=L[R_{2n-2}].
\]
By the uniqueness just proved,
\[
(R')_n=R_{2n}.
\]
Thus $R'=R$.

Finally, the assertion that a tuple is coded is uniformly $\Sigma^1_3$, while
\[
(R)_n\notin L[(R)_{n-1}]
\]
is $\Pi^1_2$. Number quantification does not raise projective complexity. Hence
\[
\forall n\in\omega\bigl((2n,(R)_n)\text{ is coded}\bigr)
\]
is $\Sigma^1_3$, and the second conjunct in $\Psi(R)$ is $\Pi^1_2$. Since $\Sigma^1_3$ is closed under conjunction with $\Pi^1_2$ formulas, $\Psi(R)$ is $\Sigma^1_3$.
\end{proof}
\subsection{The stop-and-go recursion for the second iteration}

We now work over $W[G_\omega]$. Write
\[
B^*:=B_\omega
\]
for the cumulative forbidden family obtained in the preparatory iteration. Before specifying the individual stage forcings of the main iteration, we define the main iteration together with the fixed-point classes which govern their choice. At stage $\beta$, the appropriate bookkeeping case will determine the stage forcing $\dot{\forceQ}_\beta$ and, after forcing with it, the new forbidden family $B_\beta$. 

\begin{definition}\label{MainStopAndGoRecursion}
We index the main iteration and its fixed-point phases by the nonzero countable ordinals. We recursively define
\[
(\forceP_\beta,\dot{\forceQ}_\beta\mid 1\leq\beta<\omega_1),
\]
together with the generic filters $G_\beta$, the intermediate models $M_\beta$, the derivative classes $\mathcal A_{\beta,\alpha}$, their fixed points $\mathcal A_{\beta,\infty}$, the forbidden families $B_\beta$, and countable partial well-orders $<_\beta$ of reals. Here a partial well-order means a well-order of a specified countable set of reals; $\operatorname{field}(<_\beta)$ denotes this underlying set.

\smallskip
\noindent
\emph{Initial stage.}
Let
\[
\forceP_1:=\mathbf 1,\qquad G_1:=\emptyset,\qquad M_1:=W[G_\omega],
\qquad <_1:=\emptyset.
\]
In $M_1$, carry out the fixed-$B$ construction from
Definition~\ref{AlphaAllowableAvoidingB} with $B=B^*$. For every
derivative rank $\alpha$, put
\[
\mathcal A_{1,\alpha}:=
\bigl(\mathcal A^{B^*}_\alpha\bigr)^{M_1},
\]
and let
\[
\mathcal A_{1,\infty}:=
\bigl(\mathcal A^{B^*}_\infty\bigr)^{M_1}.
\]

\smallskip
\noindent
\emph{Successor stage.}
Suppose that $1\leq\beta<\omega_1$, that $\forceP_\beta$ and $G_\beta$ have been defined, that
\[
M_\beta=W[G_\omega][G_\beta],
\]
and that $\mathcal A_{\beta,\infty}$ is the fixed-point class of phase $\beta$ in $M_\beta$. The appropriate case of the main construction described below supplies a $\forceP_\beta$-name $\dot{\forceQ}_\beta$ such that
\[
\dot{\forceQ}_\beta^{G_\beta}\in\mathcal A_{\beta,\infty}.
\]
We define
\[
\forceP_{\beta+1}:=\forceP_\beta*\dot{\forceQ}_\beta.
\]
Let $H_\beta$ be $\dot{\forceQ}_\beta^{G_\beta}$-generic over $M_\beta$, and put
\[
G_{\beta+1}:=G_\beta*H_\beta
\]
and
\[
M_{\beta+1}:=M_\beta[H_\beta]
=W[G_\omega][G_{\beta+1}].
\]
The same case of the construction also determines, in $M_{\beta+1}$, a set $B_\beta$ of newly forbidden tuples and the next partial well-order $<_{\beta+1}$. At every successor stage which is not devoted to the well-order construction we set
\[
<_{\beta+1}:=<_\beta.
\]
At a well-order stage, $<_{\beta+1}$ is defined as in the corresponding case below.

If $B_\beta=\emptyset$, no restart is performed. For every derivative rank $\alpha$, including the fixed point, we carry the same defining rule into the new model and put
\[
\mathcal A_{\beta+1,\alpha}
:=\bigl(\mathcal A_{\beta,\alpha}\bigr)^{M_{\beta+1}}.
\]
In particular,
\[
\mathcal A_{\beta+1,\infty}
=\bigl(\mathcal A_{\beta,\infty}\bigr)^{M_{\beta+1}}.
\]

If $B_\beta\neq\emptyset$, let
\[
\mathcal C^\beta_\infty:=
\bigl(\mathcal A_{\beta,\infty}\bigr)^{M_{\beta+1}}.
\]
In $M_{\beta+1}$ apply the restart from
Definition~\ref{RestartThinningOut} to $\mathcal C^\beta_\infty$ with
$D=B_\beta$, and put
\[
\mathcal A_{\beta+1,\alpha}
:=
\operatorname{Rest}(\mathcal C^\beta_\infty,B_\beta)_\alpha
\qquad(\alpha\in\operatorname{Ord}),
\]
and
\[
\mathcal A_{\beta+1,\infty}
:=
\operatorname{Rest}(\mathcal C^\beta_\infty,B_\beta)_\infty.
\]

\smallskip
\noindent
\emph{Limit stage.}
Let $\lambda<\omega_1$ be a nonzero limit ordinal and suppose that the construction has been carried out below $\lambda$. We take the finite-support direct limit
\[
\forceP_\lambda:=\operatorname{dir}\lim_{1\leq\xi<\lambda}\forceP_\xi
\]
and let
\[
G_\lambda:=\bigcup_{1\leq\xi<\lambda}G_\xi.
\]
Thus
\[
M_\lambda:=W[G_\omega][G_\lambda].
\]
Since the partial well-orders are extended only by end-extension, we put
\[
<_\lambda:=\bigcup_{1\leq\xi<\lambda}<_\xi.
\]
Reinterpret the fixed-point classes of all preceding phases in $M_\lambda$ and set
\[
\mathcal A_{\lambda,0}:=
\bigcap_{1\leq\xi<\lambda}
\bigl(\mathcal A_{\xi,\infty}\bigr)^{M_\lambda}.
\]
For $\alpha>0$, apply the derivative from Definition~\ref{AlphaAllowable} to this base class in $M_\lambda$, and let
\[
\mathcal A_{\lambda,\infty}
\]
denote the resulting fixed point.
\end{definition}

For $1\leq\beta<\omega_1$, we call $\beta$ the index of the corresponding fixed-point phase. The classes $\mathcal A_{\beta,\alpha}$ are the derivative classes within phase $\beta$, and $\mathcal A_{\beta,\infty}$ is the fixed-point class of that phase. We call their members the $(\beta,\alpha)$-allowable and $(\beta,\infty)$-allowable forcings, respectively. We continue to use \emph{stage} for a stage of the outer forcing iteration. By construction, every member of $\mathcal A_{\beta,\infty}$ avoids $B^*$ and every $B_\xi$ with $1\leq\xi<\beta$.

The final forcing is the finite-support direct limit
\[
\forceP_{\omega_1}:=\operatorname{dir}\lim_{1\leq\xi<\omega_1}\forceP_\xi,
\]
and if $G_{\omega_1}$ is generic for $\forceP_{\omega_1}$, then
\[
G_{\omega_1}=\bigcup_{1\leq\xi<\omega_1}G_\xi.
\]

At a successor stage the forcing $\dot{\forceQ}_\beta^{G_\beta}$ need not be a single coding forcing. It may itself be a countable finite-support allowable iteration. We use this entire allowable forcing as the single stage forcing at stage $\beta$ of the main iteration; its internal coding factors are not identified with separate stages of $(\forceP_\xi\mid\xi<\omega_1)$.

Each stage forcing is ccc, since it belongs to the fixed-point class of the current phase. Hence the finite-support iteration $(\forceP_\beta\mid\beta\leq\omega_1)$ is ccc. Moreover, every stage forcing has size at most $\aleph_1$, and $W[G_\omega]$ satisfies $\CH$. Inductively, each $\forceP_\beta$ has size at most $\aleph_1$, and the same is true of the finite-support limit $\forceP_{\omega_1}$. Thus $\forceP_{\omega_1}$ adds at most $\aleph_1$ reals. Since it is ccc and therefore preserves $\omega_1$, we obtain
\[
W[G_\omega][G_{\omega_1}]\models\CH.
\]

The recursion is well-defined. At a successor stage at which a new forbidden family is introduced, this follows from Lemma~\ref{RestartThinningProperties}; if no new family is introduced, the preceding defining rule is simply reinterpreted in the larger model. At a limit stage the trivial forcing belongs to every class in the intersection defining $\mathcal A_{\lambda,0}$. Moreover, if two members of this intersection use disjoint coding areas, the canonical forcing obtained by concatenating their bookkeeping functions belongs to the fixed-point class of every preceding phase, and hence to the intersection. The derivative induction therefore gives a nonempty fixed point with the product property. In particular, whenever $1\leq\beta<\gamma\leq\delta$, the classes, when computed in the common model $M_\delta$, satisfy
\[
\bigl(\mathcal A_{\gamma,\infty}\bigr)^{M_\delta}
\subseteq
\bigl(\mathcal A_{\beta,\infty}\bigr)^{M_\delta}.
\]
This is the monotonicity used below when a persistence decision made in one phase is carried through all later restarts.

We fix in $L$ a bookkeeping function $F$ of domain $\omega_1$ whose values are names for elements of $H(\omega_1)$ and such that every relevant name occurs uncountably often. Let $G_\beta$ be $\forceP_\beta$-generic and work in $W[G_\omega][G_\beta]$. We distinguish four cases.

\subsection{Continuing the diagonalization}
We reserve one bookkeeping tag for the diagonalization requirements from the preparatory construction. Suppose that at stage $\beta$ the bookkeeping presents $(m,\gamma,\dot y)$, where $\gamma<\beta$, $m$ is the G\"odel number of a $\Pi^1_3$-formula in two variables, and $\dot y$ is a $\forceP_\gamma$-name for a real. Put $y=\dot y^{G_\gamma}$.

If case~1 did not apply at the first treatment stage $\beta_m$, we use the trivial forcing. Otherwise let $x_m$ be the distinguished real associated with $m$ and let $y_m$ be the corresponding reserved real. If $y\neq y_m$, use
\[
\dot{\forceQ}_\beta^{G_\beta}=\operatorname{Code}(0,0,x_m,y,\eta),
\]
where $\eta$ is the least free coding area; if $y=y_m$, use the trivial forcing. In all cases set
\[
B_\beta=\emptyset,\qquad <_{\beta+1}=<_\beta.
\]
No restart is performed. The nontrivial factor above belongs to the fixed-point class $\mathcal A_{\beta,\infty}$ of phase $\beta$: clause~(2) of Definition~\ref{0-allowable forcing} is unchanged by the derivative, and for $y\neq y_m$ the tuple $(0,0,x_m,y)$ is not forbidden by $B_\omega$. The bookkeeping presents every relevant name cofinally often, so for every distinguished real $x_m$ it carries out case~3 throughout the main iteration. 

\subsection{Forcing \texorpdfstring{$\Pi^1_3$}{Pi13}-reduction}
We first assume that $F(\beta)=(\beta_1,\beta_2)$, where $\beta_2<\beta$, and that the $\beta_1$-th (in some previously fixed well-order $<$ of $H(\omega_2)$) $\forceP_{\beta_2}$-name of a triple of the form $(\dot n,\dot l,\dot a)$, where $\dot a$ is a nice $\forceP_{\beta_2}$-name of a real, and $\dot n,\dot l$ are nice $\forceP_{\beta_2}$-names of natural numbers, is the triple $(\dot m,\dot k,\dot x)$. We let $m=\dot m^{G_\beta}$, $k=\dot k^{G_\beta}$ and assume that $m,k$ are both G\"odel numbers of $\Pi^1_3$-formulas and $x=\dot x^{G_\beta}$ is a real. In this situation we work towards $\Pi^1_3$-reduction. We work over $W[G_\omega][G_\beta]$ as the ground model.
This model has the sequence
\[
(B_\eta\mid\eta<\beta)
\]
and the fixed-point class $\mathcal A_{\beta,\infty}$ of phase $\beta$, whose members are the $(\beta,\infty)$-allowable forcings avoiding this sequence. We distinguish three cases.

\subsubsection{Forcing \texorpdfstring{$\Pi^1_3$}{Pi13}-reduction, case 1}
We assume that in the universe $W[G_\omega][G_\beta]$ it is true that
\begin{align*}
W[G_\omega][G_\beta]\models&\forall\forceQ\bigl(\forceQ\text{ is }(\beta,\infty)\text{-allowable avoiding }\\
&(B_i\mid i<\beta)\rightarrow\forceQ\Vdash x\in A_m\bigr),
\end{align*}
then force with
\[
\dot{\forceQ}_\beta^{G_\beta}:=\operatorname{Code}(m,k,1,x,\eta).
\]
Note that this has the direct consequence that if we restrict ourselves from now on to forcings $\forceQ\in W[G_\omega][G_{\beta+1}]$ which are $(\zeta_0,\infty)$-allowable and avoid $(B_i\mid i<\zeta_0)$ for some later phase $\zeta_0>\beta$, then $x$ will remain an element of $A_m$. In particular, the pathological situation that $x\notin A_m$, $x\in A_k$ while $x$ is coded into $\vec S$ is ruled out for $(m,k,x)$.

\subsubsection{Forcing \texorpdfstring{$\Pi^1_3$}{Pi13}-reduction, case 2}
We assume that case~1 does not hold; however, we assume that in the universe $W[G_\omega][G_\beta]$,
\begin{align*}
W[G_\omega][G_\beta]\models&\forall\forceQ\bigl(\forceQ\text{ is }(\beta,\infty)\text{-allowable avoiding }\\
&(B_i\mid i<\beta)\rightarrow\forceQ\Vdash x\in A_k\bigr),
\end{align*}
then force with
\[
\dot{\forceQ}_\beta^{G_\beta}:=\operatorname{Code}(m,k,0,x,\eta).
\]

\subsubsection{Forcing \texorpdfstring{$\Pi^1_3$}{Pi13}-reduction, case 3}
In the final case we assume that neither case~1 nor case~2 applies. We obtain that there is a forcing $\forceQ$ with
\begin{align*}
W[G_\omega][G_\beta]\models\forceQ\text{ is }&(\beta,\infty)\text{-allowable avoiding $(B_i\mid i<\beta)$ and }\\
&\forceQ\Vdash x\notin A_m.
\end{align*}
Likewise we also obtain that there is a $\forceR$ such that
\begin{align*}
W[G_\omega][G_\beta]\models\forceR\text{ is }&(\beta,\infty)\text{-allowable avoiding $(B_i\mid i<\beta)$ and }\\
&\forceR\Vdash x\notin A_k.
\end{align*}
In particular there are $(\beta,\infty)$-allowable forcings avoiding $(B_i\mid i<\beta)$ which kick $x$ out of $A_m$ and $A_k$, respectively. We do want to force $x$ out of $A_m$ and $A_k$ but have to be a bit careful.

Before defining the closure operation, we record a standard observation concerning coding areas. We use it in the generic extensions under consideration in every fixed-point phase of the stop-and-go recursion.

\medskip

\noindent
\emph{Fresh coding areas observation.}
Suppose that, in the current model, $\forceP$ belongs to the fixed-point class $\mathcal A_{\beta,\infty}$ of phase $\beta$, avoids $(B_\eta\mid\eta<\beta)$, and
\[
\forceP\Vdash\varphi(a),
\]
where $a$ is a real in the current model. If $E\subseteq\omega_1$ is countable, then there is a forcing $\forceP'\in\mathcal A_{\beta,\infty}$ avoiding $(B_\eta\mid\eta<\beta)$ such that
\[
C^{\forceP'}\cap E=\emptyset
\qquad\text{and}\qquad
\forceP'\Vdash\varphi(a).
\]

Write $g^1=\langle g^1_\xi\mid\xi<\omega_1\rangle$ and, for $\rho<\omega_1$, let
\[
\mathbb C_{\geq\rho}
:=
\left(\prod_{\rho\leq\xi<\omega_1}\mathbb C(\omega_1)\right)^L.
\]
Choose $\theta<\omega_1$ such that the coordinates needed to interpret $a$, $E$, $(B_\eta\mid\eta<\beta)$, and the parameters determining phase $\beta$ and its fixed-point class $\mathcal A_{\beta,\infty}$ all lie below $\theta$. By the factorization of the Cohen product used above, the current model can be written as
\[
N[\langle g^1_\xi\mid\theta\leq\xi<\omega_1\rangle]
\]
for a suitable intermediate model $N$ containing these parameters. Hence there are a $\mathbb C_{\geq\theta}$-name $\dot{\forceP}$ and a condition $p\in\mathbb C_{\geq\theta}$ such that
\begin{equation*}
\tag{1}
p\Vdash_{\mathbb C_{\geq\theta}}
\left[
\begin{array}{l}
\dot{\forceP}\in\mathcal A_{\beta,\infty},\\[1mm]
\dot{\forceP}\text{ avoids }(B_\eta\mid\eta<\beta),\\[1mm]
\dot{\forceP}\Vdash\varphi(a)
\end{array}
\right].
\end{equation*}

Choose
\[
\zeta>\max\{\theta,\sup E\}
\]
and let
\[
\pi_{\theta,\zeta}:\mathbb C_{\geq\theta}\cong\mathbb C_{\geq\zeta}
\]
be a coordinate isomorphism. Transporting $p$, $\dot{\forceP}$, and the coding-area coordinates occurring in $\dot{\forceP}$ along $\pi_{\theta,\zeta}$ gives a condition $p_\zeta$ and a name $\dot{\forceP}_\zeta$ such that
\begin{equation*}
\tag{2}
p_\zeta\Vdash_{\mathbb C_{\geq\zeta}}
\left[
\begin{array}{l}
\dot{\forceP}_\zeta\in\mathcal A_{\beta,\infty},\\[1mm]
\dot{\forceP}_\zeta\text{ avoids }(B_\eta\mid\eta<\beta),\\[1mm]
C^{\dot{\forceP}_\zeta}\subseteq[\zeta,\omega_1),\\[1mm]
\dot{\forceP}_\zeta\Vdash\varphi(a)
\end{array}
\right].
\end{equation*}
The definitions of membership in the fixed-point class $\mathcal A_{\beta,\infty}$ of phase $\beta$ and of avoidance are invariant under this coordinate shift, while the parameters in $N$ remain fixed. Since $\mathbb C_{\geq\zeta}$ is weakly homogeneous, (2) implies
\begin{equation*}
\tag{3}
\mathbf 1_{\mathbb C_{\geq\zeta}}\Vdash
\exists\forceQ\left[
\begin{array}{l}
\forceQ\in\mathcal A_{\beta,\infty},\\[1mm]
\forceQ\text{ avoids }(B_\eta\mid\eta<\beta),\\[1mm]
C^{\forceQ}\subseteq[\zeta,\omega_1),\\[1mm]
\forceQ\Vdash\varphi(a)
\end{array}
\right].
\end{equation*}
Evaluating (3) in the actual tail generic yields the required $\forceP'$. Indeed,
\[
E\subseteq\zeta
\qquad\text{and}\qquad
C^{\forceP'}\subseteq[\zeta,\omega_1),
\]
so $C^{\forceP'}\cap E=\emptyset$.

We now define the closure operation. Let
\[
\forceR=((\forceR_\eta,\dot{\forceQ}_\eta)\mid\eta<\gamma)
\]
be a $(\beta,\infty)$-allowable forcing avoiding $(B_\eta\mid\eta<\beta)$. Suppose that $\nu<\gamma$ is a stage at which the bookkeeping associated with $\forceR$ considers a real $a$ and two $\Pi^1_3$-sets $A_m,A_k$. Work in the corresponding $\forceR_\nu$-generic extension and write
\[
\forceQ_\nu=\dot{\forceQ}_\nu^{G_\nu}.
\]

Assume that there are $(\beta,\infty)$-allowable forcings $\forceP_m$ and $\forceP_k$, both avoiding $(B_\eta\mid\eta<\beta)$, such that
\[
\forceP_m\Vdash a\notin A_m
\qquad\text{and}\qquad
\forceP_k\Vdash a\notin A_k.
\]
Using the fresh coding areas observation twice, we may assume that
\[
C^{\forceQ_\nu},\qquad C^{\forceP_m},\qquad C^{\forceP_k}
\]
are pairwise disjoint. By the product property of the fixed-point class $\mathcal A_{\beta,\infty}$ of phase $\beta$, the product $\forceP_m\times\forceP_k$ densely embeds into a $(\beta,\infty)$-allowable forcing $\forceR'_\nu$ avoiding $(B_\eta\mid\eta<\beta)$. Since the assertions
\[
a\notin A_m
\qquad\text{and}\qquad
a\notin A_k
\]
are $\Sigma^1_3$, they are upward absolute to the product extension. Hence
\[
\forceR'_\nu\Vdash a\notin A_m\cup A_k.
\]
The forcing $\forceR'_\nu$ may be obtained by concatenating the bookkeeping functions for $\forceP_m$ and $\forceP_k$, and therefore uses no coding areas other than those used by these two forcings. In particular,
\[
C^{\forceR'_\nu}\cap C^{\forceQ_\nu}=\emptyset.
\]

Applying the product property once more, choose a $(\beta,\infty)$-allowable forcing $C^1(\forceQ_\nu)$, avoiding $(B_\eta\mid\eta<\beta)$, into which
\[
\forceQ_\nu\times\forceR'_\nu
\]
densely embeds. If no such pair $\forceP_m,\forceP_k$ exists, let
\[
C^1(\forceQ_\nu):=\forceQ_\nu.
\]

Applying this construction at every stage $\nu<\gamma$ and taking the finite-support iteration with these replacement factors gives the first closure
\[
C^1(\forceR):=
\Asterisk_{\nu<\gamma}C^1(\dot{\forceQ}_\nu).
\]
Recursively let
\[
C^{n+1}(\forceR):=C^1(C^n(\forceR))
\]
and let $C^\omega(\forceR)$ be the finite-support direct limit of $(C^n(\forceR)\mid n<\omega)$. We call $C^\omega(\forceR)$ the closure of $\forceR$.

\begin{lemma}\label{ClosureLemma}
Let $\forceR$ be a $(\beta,\infty)$-allowable forcing avoiding $(B_\eta\mid\eta<\beta)$, relative to a bookkeeping function $F'$. Then $C^\omega(\forceR)$ is a $(\beta,\infty)$-allowable forcing avoiding $(B_\eta\mid\eta<\beta)$, relative to a suitable bookkeeping function $F''$.

Moreover, whenever $F''$ considers a triple $(m,k,a)$ and there are $(\beta,\infty)$-allowable forcings avoiding $(B_\eta\mid\eta<\beta)$ which force $a\notin A_m$ and $a\notin A_k$, respectively, then
\[
C^\omega(\forceR)\Vdash a\notin A_m\cup A_k.
\]
\end{lemma}

\begin{proof}
At every stage of the closure construction, the fresh coding areas observation allows us to choose the auxiliary forcings with coding areas disjoint from the current factor. The product property therefore provides a replacement factor with the same degree of allowability. The replacement also avoids $(B_\eta\mid\eta<\beta)$, since all its coding factors already occur in forcings which avoid this family.

It follows inductively that every $C^n(\forceR)$ is $(\beta,\infty)$-allowable and avoids $(B_\eta\mid\eta<\beta)$. The same is true of the finite-support direct limit $C^\omega(\forceR)$.

Finally, every triple appearing in the bookkeeping of $C^\omega(\forceR)$ already appears at some finite stage $C^n(\forceR)$. If the corresponding real can be forced out of both $\Pi^1_3$-sets, the construction of $C^{n+1}(\forceR)$ inserts a forcing which forces it out of their union. This forcing remains a complete part of every later closure and hence of $C^\omega(\forceR)$.
\end{proof}
Returning to our iteration at stage $\beta$, let $C^\omega(\forceQ)$ and $C^\omega(\forceR)$ be the closures of the $<$-least witnesses described above. Using the fresh coding areas observation throughout the two closure constructions, we may arrange that
\[
C^{C^\omega(\forceQ)}\cap C^{C^\omega(\forceR)}=\emptyset.
\]
The product property now gives a $(\beta,\infty)$-allowable forcing avoiding $(B_i\mid i<\beta)$ into which
\[
C^\omega(\forceQ)\times C^\omega(\forceR)
\]
densely embeds. We use this forcing as $\dot{\forceQ}_\beta^{G_\beta}$. Since the two factors force $x\notin A_m$ and $x\notin A_k$, respectively, upward absoluteness of these $\Sigma^1_3$ statements implies that
\[
\dot{\forceQ}_\beta^{G_\beta}\Vdash x\notin A_m\cup A_k.
\]

For all three cases, let $B_\beta:=\emptyset$. No restart is performed and, in $M_{\beta+1}$, we put
\[
\mathcal A_{\beta+1,\infty}
:=\bigl(\mathcal A_{\beta,\infty}\bigr)^{M_{\beta+1}}.
\]

This ends the definition of the iteration in this case and we shall show that, if $G_{\omega_1}$ denotes a generic filter for the forcing $\forceP_{\omega_1}$, which is defined as the direct limit of the forcings $\forceP_\beta$, then the resulting universe $W[G_\omega][G_{\omega_1}]$ satisfies the $\Pi^1_3$-reduction property.
\subsection{Forcing \texorpdfstring{$\Sigma^1_4$}{Sigma14}-uniformization}

This section deals with forcing the $\Sigma^1_4$-uniformization property. First we note that we need to define how to force $\Sigma^1_4$-uniformization in a different way to \cite{BPFA and uniformization}. This is necessary, as the latter machinery can not be combined with forcing $\Pi^1_3$-reduction. Indeed both methods exclude themselves and so something new is in demand.
For every integer $n\in\omega$ we define the set
\[
f_n:=\{(x,y)\mid\exists a_0\forall a_1((n,x,y,a_0,a_1)\text{ is not coded into $\vec S$})\}
\]
and will work towards a universe where for every $n\in\omega$, if $n$ is the G\"odel number of a $\Sigma^1_4$-set $A_n$ in the plane then $f_n$ uniformizes $A_n$. Note that $f_n$ is a $\Sigma^1_4$-formula, hence $\Sigma^1_4$-uniformization would hold.

Assume that $F(\beta)$ hands us a G\"odel number $m$ of a $\Sigma^1_4$-formula
\[
\varphi_m=\exists a_0\forall a_1\,\psi_m(x,y,a_0,a_1),
\]
where $\psi_m$ is $\Sigma^1_2$, and reals $x,y',a'_0$. We work in $W[G_\omega][G_\beta]$ with the fixed-point class $\mathcal A_{\beta,\infty}$ of phase $\beta$. Fix once and for all
\[
a_1^*:=0^\omega.
\]

We first deal with the situation in which a value has already been selected. Suppose that at an earlier $\Sigma^1_4$-stage we declared
\[
f_m(x)=y.
\]
If $y'\neq y$, let $\eta<\omega_1$ be the least coding area which is free and put
\[
\dot{\forceQ}_\beta^{G_\beta}
:=\operatorname{Code}(m,x,y',a'_0,a_1^*,\eta).
\]
If $y'=y$, we use the trivial forcing. In either case let $B_\beta:=\emptyset$. No restart is performed. Since every relevant quadruple $(m,x,y',a_0)$ occurs uncountably often in the bookkeeping, once $f_m(x)=y$ has been selected, every pair $y'\neq y$ and $a_0$ is eventually treated in this way. Thus in the final model
\[
\forall y'\neq y\,\forall a_0\,
\bigl((m,x,y',a_0,a_1^*)\text{ is coded into }\vec S\bigr).
\]

It remains to describe what happens when no value for $f_m(x)$ has yet been selected.

\subsubsection{Case 1}

We ask whether there are a real $y$ and a $(\beta,\infty)$-allowable forcing
\[
\forceR=((\forceR_\eta,\dot{\forceQ}_\eta)\mid\eta<\gamma)
\]
avoiding $(B_\eta\mid\eta<\beta)$ such that, if $H$ is $\forceR$-generic over $W[G_\omega][G_\beta]$, then $W[G_\omega][G_\beta][H]$ contains a real $a_0$ for which every further $(\beta,\infty)$-allowable forcing $\forceR'$ avoiding $(B_\eta\mid\eta<\beta)$ satisfies
\[
\forceR'\Vdash\forall a_1\,\psi_m(x,y,a_0,a_1).
\]
In this situation we fix the $<$-least such witness and declare
\[
f_m(x):=y.
\]
We would want to use $\forceR$ to add the corresponding real $a_0$, but we need to ensure that the forcing used at this stage does not interfere with the reduction construction.

We therefore first pass to the closure
\[
\mathbb T:=C^\omega(\forceR).
\]
The forcing $\mathbb T$ still adds the real $a_0$, and by Lemma~\ref{ClosureLemma} it has the same degree of allowability and avoids the same forbidden families as $\forceR$.

We explain why passing to the closure prevents this forcing from interfering with the proof of $\Pi^1_3$-reduction. Consider a triple $(p,q,a)$ which occurs in the reduction bookkeeping of $C^\omega(\forceR)$. If $a\notin A_p\cup A_q$ already holds when the triple is considered, then this remains true in every later forcing extension, since this is a $\Sigma^1_3$-statement, and any coding decision concerning $a$ is irrelevant to the reducing sets.

Suppose next that one of the two persistence alternatives in the definition of the fixed-point class $\mathcal A_{\beta,\infty}$ of phase $\beta$ holds. For example, suppose that every continuation from this class forces $a\in A_p$. The definition of allowability then chooses the coding direction which places $a$ into the reducing set associated with $A_p$. This decision remains correct after all later restarts of the thinning-out process, since the fixed-point class of every later phase is a subclass of $\mathcal A_{\beta,\infty}$. Thus a persistence statement which holds in phase $\beta$ continues to hold throughout all later phases.

It remains to consider the case in which neither persistence alternative holds. There are then forcings $\forceP_p$ and $\forceP_q$ in the fixed-point class of phase $\beta$, avoiding the current forbidden families, such that
\[
\forceP_p\Vdash a\notin A_p
\qquad\text{and}\qquad
\forceP_q\Vdash a\notin A_q.
\]
By the construction of the closure, after moving their coding areas to fresh coordinates and applying the product property, the next closure step inserts a forcing which forces
\[
a\notin A_p\cup A_q.
\]
This is again a $\Sigma^1_3$-statement and therefore remains true in every subsequent generic extension. Consequently, any provisional coding decision made for the triple $(p,q,a)$ becomes irrelevant to the final reducing sets.

Every triple appearing in the bookkeeping of $C^\omega(\forceR)$ occurs at some finite stage of the closure construction and is treated at the following stage. Hence every reduction decision occurring inside the forcing used to add $a_0$ is either justified by a persistence statement which survives all later thinning-out, or is rendered irrelevant by forcing the corresponding real out of both sets. It follows that using $C^\omega(\forceR)$ at the present $\Sigma^1_4$-uniformization stage does not disturb the eventual proof of $\Pi^1_3$-reduction.

We now include the real which will witness the definition of $f_m(x)$ in the same stage forcing. Fix a G\"odel number $l$ of a $\Sigma^1_4$-formula which defines the empty set, and fix reals $u,v,b_0,b_1$. Since $A_l=\emptyset$, no tuple with first coordinate $l$ is ever reserved by a successful uniformization stage. By the arity-sensitive coding convention, the other forbidden families do not interfere with this five-tuple. Let $\dot{\mathbb S}$ be the $\mathbb T$-name which, in a $\mathbb T$-generic extension, is interpreted as
\[
\operatorname{Code}(l,u,v,b_0,b_1,\eta^*),
\]
where $\eta^*$ is the least coding area which is free in that extension. We define the complete stage forcing by
\[
\dot{\forceQ}_\beta^{G_\beta}:=\mathbb T*\dot{\mathbb S}.
\]
The second factor is one of the coding factors from clause~(3) of Definition~\ref{0-allowable forcing}. The derivative leaves this clause unchanged, so appending this factor at a free coding area preserves $(\beta,\infty)$-allowability and avoidance of the current forbidden families. Let
\[
c_0(x,y,m)
\]
be the almost-disjoint generic real added by $\dot{\mathbb S}$. Thus $c_0(x,y,m)\in M_{\beta+1}$.

In $M_{\beta+1}$ define
\[
B_\beta:=\{(m,x,y,c_0(x,y,m),\dot z)\mid
\dot z\text{ is a nice name for a real in an allowable continuation}\}.
\]
We define $\mathcal A_{\beta+1,\infty}$ by the successor clause of Definition~\ref{MainStopAndGoRecursion}: restrict the fixed-point class of phase $\beta$ to forcings which avoid $B_\beta$ and apply the derivative until it reaches its fixed point. Hence no later forcing uses
\[
\operatorname{Code}(m,x,y,c_0(x,y,m),z,\eta)
\]
for any real $z$ and any coding area $\eta$. Since $c_0(x,y,m)$ is new at this stage, no such tuple was coded earlier either. Therefore the final model satisfies
\[
\forall a_1\,
\bigl((m,x,y,c_0(x,y,m),a_1)\text{ is not coded into }\vec S\bigr).
\]
Thus
\[
\exists a_0\forall a_1\,
\bigl((m,x,y,a_0,a_1)\text{ is not coded into }\vec S\bigr),
\]
with $c_0(x,y,m)$ as a witness for $a_0$.

At every later $\Sigma^1_4$-stage at which the bookkeeping presents $(m,x,y',a_0)$ with $y'\neq y$, the preliminary clause above uses
\[
\operatorname{Code}(m,x,y',a_0,a_1^*,\eta)
\]
for a fresh coding area. Consequently, in the final model every $y'\neq y$ satisfies
\[
\forall a_0\exists a_1\,
\bigl((m,x,y',a_0,a_1)\text{ is coded into }\vec S\bigr),
\]
and hence fails the formula defining $f_m(x)$. Thus $y$ is the unique real satisfying that formula.

\subsubsection{Case 2}

Assume that case~1 fails. The negation of case~1 says that no candidate remains persistent with respect to the fixed-point class $\mathcal A_{\beta,\infty}$ of phase $\beta$. In particular, apply this to the bookkeeping values $y'$ and $a'_0$, viewed as already present in $W[G_\omega][G_\beta]$. There is then a $(\beta,\infty)$-allowable forcing $\forceQ'$ avoiding $(B_\eta\mid\eta<\beta)$ such that
\[
\forceQ'\Vdash\exists a_1\neg\psi_m(x,y',a'_0,a_1).
\]
Choose the $<$-least such forcing and put
\[
\dot{\forceQ}_\beta^{G_\beta}:=C^\omega(\forceQ').
\]
The closure remains $(\beta,\infty)$-allowable and avoids the current forbidden families. The construction of the closure gives a canonical complete embedding of $\forceQ'$ into $C^\omega(\forceQ')$, so it also forces
\[
\exists a_1\neg\psi_m(x,y',a'_0,a_1).
\]
This is a $\Sigma^1_3$-statement and hence remains true in every later generic extension.

We let $B_\beta:=\emptyset$. By Definition~\ref{MainStopAndGoRecursion}, no restart is performed and
\[
\mathcal A_{\beta+1,\infty}
:=\bigl(\mathcal A_{\beta,\infty}\bigr)^{M_{\beta+1}}.
\]

This ends the definition of the stages devoted to $\Sigma^1_4$-uniformization.

\subsection{Forcing a good \texorpdfstring{$\Sigma^1_5$}{Sigma15}-well-order}
Our methods allow for an additional layer of complexity which we can use to force a good $\Sigma^1_5$-well-order of the reals.
We single out coding forcings of the form
\[
\operatorname{Code}(1,1,x,z,\eta)
\]
to be our tool which eventually should yield the good $\Sigma^1_5$-well-order.
The idea is to use these forcings in such a way that
\begin{align*}
(1,1,x,z)&\text{ is coded ``cofinally often''}\\
&\quad\Longleftrightarrow\quad
\text{$z$ codes the initial segment of the well-order below $x$.}
\end{align*}
Here being coded ``cofinally often'' is just an abbreviation for the formula ``$\exists r_0$ ($r_0$ codes the pair $(x,z)$ and $\forall r_1\exists r_2$ ($r_2\notin L[r_1]$ and $r_2$ witnesses that $(1,1,x,z)$ is coded))''. Note that the latter formula is of the form $\exists r_0(\Delta^1_1\land\forall r_1\exists r_2(\Pi^1_2\rightarrow\Pi^1_2))$, which is $\Sigma^1_5$.
Thus the well-order we shall define will become a good $\Sigma^1_5$-well-order.

At stage $\beta$ we work in $W[G_\omega][G_\beta]$ with the fixed-point class $\mathcal A_{\beta,\infty}$ of phase $\beta$ and the current countable partial well-order $<_\beta$. We arrange the bookkeeping so that every pair of names for reals occurs uncountably often in the present case. Let the bookkeeping at $\beta$ hand us $(\dot m,\dot x,\dot z)$, let
\[
x=\dot x^{G_\beta},\qquad z=\dot z^{G_\beta},\qquad m=\dot m^{G_\beta},
\]
and assume that $m$ is the G\"odel number of a $\Sigma^1_5$-formula.

We first define $<_{\beta+1}$. If $x\in\operatorname{field}(<_\beta)$, set
\[
<_{\beta+1}:=<_\beta.
\]
If $x\notin\operatorname{field}(<_\beta)$, adjoin $x$ as a new largest element. Thus the field of $<_{\beta+1}$ is
\[
\operatorname{field}(<_\beta)\cup\{x\},
\]
and
\[
<_{\beta+1}
:=<_\beta\cup
\{(u,x):u\in\operatorname{field}(<_\beta)\}.
\]
In particular $<_{\beta+1}$ is an end-extension of $<_\beta$, so predecessor sets of reals already in the field never change later in the construction.

We now inspect the real $z$ supplied by the bookkeeping. If
\[
\{(z)_n:n\in\omega\}
=
\{u:u<_{\beta+1}x\},
\]
we force with
\[
\dot{\forceQ}_\beta^{G_\beta}
:=\operatorname{Code}(1,1,x,z,\eta)
\]
for the least coding area $\eta$ which is free. If $z$ does not code this predecessor set, we let $\dot{\forceQ}_\beta^{G_\beta}$ be the trivial forcing. Since every pair $(x,z)$ is presented uncountably often, every real $z$ which correctly codes the predecessor set of $x$ is coded at cofinally many stages after $x$ has entered the order.

At the same time we exclude incorrect predecessor codes. In $M_{\beta+1}$ let
\begin{align*}
B:=\{(1,1,x',\dot y)\mid {}&x'\in\operatorname{field}(<_{\beta+1})\text{ and there is a forcing }\forceR\in
(\mathcal A_{\beta,\infty})^{M_{\beta+1}}\text{ such that}\\
&\dot y\text{ is a nice $\forceR$-name for a real and }\forceR\Vdash
\{(\dot y)_n:n\in\omega\}\neq\{u:u<_{\beta+1}x'\}\}.
\end{align*}
We set
\[
B_\beta:=B.
\]
Thus, once a real is in the field of the partial well-order, no later forcing can code a real which gives an incorrect predecessor set for it.

Define $\mathcal A_{\beta+1,\infty}$ by the successor clause of Definition~\ref{MainStopAndGoRecursion}, using $B_\beta$ as the new set of forbidden tuples.
\section{Discussion of the universe}

In this section we show that the just defined universe has the desired properties.
Let $G_{\omega_1}$ denote a generic filter over $W[G_\omega]$ for the just defined iteration. We argue in $W[G_\omega][G_{\omega_1}]$ from now on.

\begin{lemma}\label{FailureUniformizationFinalModel}
In $W[G_\omega][G_{\omega_1}]$, the $\Pi^1_3$-uniformization property fails.
\end{lemma}

\begin{proof}
Recall the set
\[
A:=\{(x,y)\in(\omega^\omega)^2:(0,0,x,y)\text{ is not coded into }\vec S\},
\]
which is $\Pi^1_3$. The preparatory iteration produced the forbidden family $B_\omega$ and, for every G\"odel number $m$ for which case~1 applied at the first treatment stage $\beta_m$, a distinguished real $x_m$.

By construction, every forcing used in the main iteration avoids $B_\omega$. Moreover, the diagonalization part of the bookkeeping revisits every distinguished real $x_m$ cofinally often and treats every relevant real $y\neq y_m$ according to case~3. Hence the hypotheses of Lemma~\ref{PropertiesOfIntermediateModel} are satisfied by the main continuation. It follows that no $\Pi^1_3$-set in the plane uniformizes $A$ in the final model. Since $A$ itself is $\Pi^1_3$, we therefore obtain
\[
W[G_\omega][G_{\omega_1}]\models\neg\Pi^1_3\text{-uniformization}.
\]
\end{proof}

\subsection{\texorpdfstring{$\Pi^1_3$}{Pi13}-reduction holds in our universe}
We first argue why $\Pi^1_3$-reduction holds in $W[G_\omega][G_{\omega_1}]$.

For every pair $(m,k)\in\omega^2$, we define in a first step
\[
D^0_{m,k}:=\{x\in\omega^\omega:(m,k,0,x)\text{ is not coded into the }\vec S\text{-sequence}\}
\]
and
\[
D^1_{m,k}:=\{x\in\omega^\omega:(m,k,1,x)\text{ is not coded into the }\vec S\text{-sequence}\}.
\]
Now both sets will not necessarily reduce $A_m$ and $A_k$, since the preparatory passage from $W$ to $W[G_\omega]$ may have produced coding witnesses which are not intended to enter the later reduction construction. We use the formula $\Psi$ from Lemma~\ref{DefinabilityR} to separate these witnesses from the later codes. Since every marker tuple $(2n,\dot z)$ relevant to $\Psi$ is forbidden after the corresponding preparatory stage, the later main iteration adds no new marker code. Hence in $W[G_\omega][G_{\omega_1}]$ the formula $\Psi$ still defines the same unique real $R$.

For $i\in\{0,1\}$ define
\[
E^i_{m,k}:=\Bigl\{x\in\omega^\omega:\forall R\Bigl(\Psi(R)\rightarrow
\neg\exists r\bigl(r\notin L[R]\land r\text{ witnesses that }(m,k,i,x)
\text{ is coded into }\vec S\bigr)\Bigr)\Bigr\}.
\]
Both $E^0_{m,k}$ and $E^1_{m,k}$ are $\Pi^1_3$-definable. Indeed, $\Psi(R)$ is $\Sigma^1_3$ by Lemma~\ref{DefinabilityR}, while
\[
r\notin L[R]\land r\text{ witnesses that }(m,k,i,x)\text{ is coded into }\vec S
\]
is $\Pi^1_2$. Hence the existential quantification over $r$ gives a $\Sigma^1_3$ formula, and therefore
\[
\Psi(R)\rightarrow
\neg\exists r\bigl(r\notin L[R]\land r\text{ witnesses that }(m,k,i,x)
\text{ is coded into }\vec S\bigr)
\]
is $\Pi^1_3$. The outer universal quantifier over $R$ can be absorbed into the leading universal real quantifier block, so
\[
\forall R\Bigl(\Psi(R)\rightarrow
\neg\exists r\bigl(r\notin L[R]\land r\text{ witnesses that }(m,k,i,x)
\text{ is coded into }\vec S\bigr)\Bigr)
\]
is $\Pi^1_3$.

Our goal is to show that for every pair $(m,k)$ the sets $E^0_{m,k}\cap A_m$ and $E^1_{m,k}\cap A_k$ reduce the pair of $\Pi^1_3$-sets $A_m$ and $A_k$.

\begin{lemma}
In $W[G_\omega][G_{\omega_1}]$, for every pair $(m,k)$, $m,k\in\omega$, and corresponding $\Pi^1_3$-sets $A_m$ and $A_k$:
\begin{enumerate}
\item[(a)] $E^0_{m,k}\cap A_m$ and $E^1_{m,k}\cap A_k$ are disjoint.
\item[(b)] $(E^0_{m,k}\cap A_m)\cup(E^1_{m,k}\cap A_k)=A_m\cup A_k$.
\item[(c)] $E^0_{m,k}\cap A_m$ and $E^1_{m,k}\cap A_k$ are $\Pi^1_3$-definable.
\end{enumerate}
\end{lemma}

\begin{proof}
We prove~(a) first. We argue in $W[G_\omega][G_{\omega_1}]$ but ignore the codes created in $W[G_\omega]$. This is justified by the above discussion. If $x$ is an arbitrary real in $A_m\cap A_k$, there will be a least stage $\beta$ above $\omega$ such that $F$ at stage $\beta$ considers a triple of names which itself corresponds to the triple $(m,k,x)$. As $x\in A_m\cap A_k$, we know that case~1 or case~2 in the definition of our iteration in the case where we force towards $\Pi^1_3$-reduction must have applied. We argue for case~1, as case~2 is similar. In case~1 the application of $\operatorname{Code}(m,k,1,x,\eta)$ codes $(m,k,1,x)$ into $\vec S$, while ensuring that for all future extensions, $x$ will remain an element of $A_m$. The rules of the iteration also tell us that $(m,k,0,x)$ will never be coded into $\vec S$ by a later factor of the iteration. Thus $x\in E^0_{m,k}\cap A_m$. Since $(m,k,1,x)$ is coded, $x\notin E^1_{m,k}$, and therefore $E^0_{m,k}\cap A_m$ and $E^1_{m,k}\cap A_k$ are disjoint.

To prove~(b), let $x$ be an arbitrary element of $A_m\cup A_k$. Let $\beta>\omega$ be the stage of the iteration beyond $W[G_\omega]$ where the triple $(m,k,x)$ is considered first. As $x\in A_m\cup A_k$, either case~1 or case~2 in the definition of the forcing for $\Pi^1_3$-reduction was applied at stage $\beta$.

Assume first that it was case~1. Then, as argued above, $x\in A_m$ will remain true for the rest of the iteration, case~1 codes $(m,k,1,x)$, and we will never code $(m,k,0,x)$ into $\vec S$ at a later stage of our iteration. Hence $x\in A_m\cap E^0_{m,k}$. If at stage $\beta$ case~2 applied, then case~2 codes $(m,k,0,x)$, while $(m,k,1,x)$ is never coded later; hence $x\in E^1_{m,k}\cap A_k$. Thus either $x\in E^0_{m,k}\cap A_m$ or $x\in E^1_{m,k}\cap A_k$, and we are finished.

Clause~(c) holds by the calculation above.
\end{proof}

\subsection{Proof of \texorpdfstring{$\Sigma^1_4$}{Sigma14}-uniformization}

We claim that adding the two cases to our definition of the iteration will force $\Sigma^1_4$-uniformization in the final model.

\begin{theorem}
Let $G_{\omega_1}$ be a $\forceP_{\omega_1}$-generic filter, where $\forceP_{\omega_1}$ denotes the limit of the iteration $(\forceP_\alpha,\forceQ_\beta\mid\beta<\omega_1,\alpha\leq\omega_1)$.
Then in $W[G_\omega][G_{\omega_1}]$ the $\Sigma^1_4$-uniformization property holds, the $\Pi^1_3$-reduction property holds and the $\Pi^1_3$-uniformization property fails.
\end{theorem}

\begin{proof}
We first prove the $\Sigma^1_4$-uniformization property. Let $m$ be the G\"odel number of a $\Sigma^1_4$-formula in two free variables
\[
\varphi_m(x,y)\equiv\exists a_0\forall a_1\,\psi_m(x,y,a_0,a_1),
\]
where $\psi_m$ is $\Sigma^1_2$, and let $x$ be an arbitrary real in the final model.

Assume first that no case~1 stage ever selects a value for $f_m(x)$. Then case~1 fails whenever the bookkeeping treats $(m,x,y,a_0)$. By case~2, for every pair of reals $y,a_0$ which appears in the final model, a later stage forces
\[
\exists a_1\neg\psi_m(x,y,a_0,a_1).
\]
Every relevant name occurs uncountably often, so every such pair is eventually treated after it appears. Since the displayed assertion is $\Sigma^1_3$, it remains true in the final model. Hence
\[
\forall y\,\forall a_0\,\exists a_1\neg\psi_m(x,y,a_0,a_1),
\]
and the $x$-section of $A_m$ is empty.

Suppose instead that at some stage $\beta$ case~1 selects a value
\[
f_m(x)=y.
\]
At that stage the forcing adds a real $a_0$ such that every later allowable continuation in the fixed-point classes of subsequent phases preserves
\[
\forall a_1\,\psi_m(x,y,a_0,a_1).
\]
Since the fixed-point class of every later phase is a subclass of the fixed-point class of phase $\beta$, this remains true in the final model. Thus
\[
(x,y)\in A_m.
\]
The same stage also adds the real $c_0(x,y,m)$ and then places every tuple
\[
(m,x,y,c_0(x,y,m),\dot z)
\]
in the forbidden family. Consequently,
\[
\forall a_1\,
\bigl((m,x,y,c_0(x,y,m),a_1)\text{ is not coded into }\vec S\bigr),
\]
so
\[
\exists a_0\forall a_1\,
\bigl((m,x,y,a_0,a_1)\text{ is not coded into }\vec S\bigr).
\]

On the other hand, after $y$ has been selected, the bookkeeping treats every quadruple $(m,x,y',a_0)$ with $y'\neq y$ and uses
\[
\operatorname{Code}(m,x,y',a_0,a_1^*,\eta)
\]
at a fresh coding area. Hence, in the final model, for every $y'\neq y$,
\[
\forall a_0\exists a_1\,
\bigl((m,x,y',a_0,a_1)\text{ is coded into }\vec S\bigr).
\]
Therefore $y$ is the unique real satisfying the $\Sigma^1_4$-formula defining $f_m(x)$. Since $m$ and $x$ were arbitrary, $\Sigma^1_4$-uniformization holds.

The $\Pi^1_3$-reduction property follows from the preceding reduction lemma, and the failure of $\Pi^1_3$-uniformization follows from Lemma~\ref{FailureUniformizationFinalModel}. Thus all three assertions of the theorem hold.
\end{proof}

\subsection{Proof of the good \texorpdfstring{$\Sigma^1_5$}{Sigma15}-wellorder}

\begin{lemma}
Let $W[G_\omega][G_{\omega_1}]$ be the result if we run our definition of the iteration for $\omega_1$-many steps over $W[G_\omega]$.
Then the reals of $W[G_\omega][G_{\omega_1}]$ have a good $\Sigma^1_5$-well-order of the reals. Hence $\Sigma^1_m$-uniformization holds for every $m\geq5$.
\end{lemma}
\begin{proof}
In $W[G_\omega][G_{\omega_1}]$ define
\[
<\;:=\bigcup_{1\leq\beta<\omega_1}<_\beta.
\]
Every $<_{\beta+1}$ is an end-extension of $<_\beta$, and at limit stages we take unions. Hence $<$ is a well-order and every proper initial segment of $<$ is countable. Every real in the final model is the interpretation of a name appearing at some countable stage of the ccc iteration, and the bookkeeping presents that name at a later well-order stage. Thus every real eventually belongs to the field of $<$. By the calculation following Definition~\ref{MainStopAndGoRecursion}, the final model satisfies $\CH$, so there are exactly $\omega_1$ reals. It follows that $<$ has order type $\omega_1$.

Now let $x,z$ be reals and suppose that $z$ codes the set of $<$-predecessors of $x$. Once $x$ has entered the field of the partial well-order, this predecessor set never changes. Since the bookkeeping presents the pair $(x,z)$ uncountably often, the forcing $\operatorname{Code}(1,1,x,z,\eta)$ is used cofinally often. Hence $(1,1,x,z)$ satisfies the $\Sigma^1_5$-formula ``$(1,1,x,z)$ is coded cofinally often''.

Conversely, suppose that $z$ does not code the $<$-initial segment below $x$. There is then a stage $\beta<\omega_1$ such that $x\in\operatorname{field}(<_{\beta+1})$ and
\[
\{(z)_n:n\in\omega\}\neq\{u:u<_{\beta+1}x\}.
\]
At that stage the tuple $(1,1,x,z)$ belongs to $B_\beta$, and the fixed-point class of every later phase avoids it. Thus no later stage codes this tuple. The part of the iteration below $\beta+1$ is countable, so its coding reals can be coded by a single real $r$. Consequently every witness that $(1,1,x,z)$ was coded before stage $\beta+1$ belongs to $L[r]$, while there are no later witnesses. Hence there is no $s\notin L[r]$ which witnesses that $(1,1,x,z)$ is coded into $\vec S$. Thus $(1,1,x,z)$ is not coded cofinally often.

We have therefore shown that
\[
(1,1,x,z)\text{ is coded cofinally often}
\quad\Longleftrightarrow\quad
z\text{ codes the set of $<$-predecessors of }x.
\]
This is a $\Sigma^1_5$ definition of the initial segments of $<$, so $<$ is a good $\Sigma^1_5$-well-order.
\end{proof}

\section{The second main theorem}
In this section we want to sketch a proof of the second main theorem. Its proof can be seen as an easier variant of the proof from the first main theorem.
Recall what we are aiming for:
\begin{theorem}
There is a generic extension of $L$ where the
$\Pi^1_3$-uniformization property holds and where $\Sigma^1_n$-uniformization holds for $n \ge 4$.
\end{theorem}
The theorem is proved using an $\omega_1$-length iteration of allowable forcings; the notion of allowable needs to altered though for our new task. It should accomodate the following tasks:
\begin{enumerate}
\item Forcing $\Pi^1_3$-uniformization

\item Forcing a good $\Sigma^1_5$ well-order of the reals.
\end{enumerate}

\subsection{Allowable forcings for the construction}

Hence we define a new variant of 0-allowable forcings designed to carry out such a proof.

\begin{definition}
    
    Let $W$ be our ground model. Let $\alpha < \omega_1$ and let $F\in W$, $F: \alpha \rightarrow W$ be a bookkeeping function.
    
    A finite support iteration $\forceP=(\forceP_{\beta}\,:\, {\beta< \alpha})$ is called allowable (relative to the bookkeeping function $F$) if the function $F: \alpha \rightarrow W$ determines $\forceP$ inductively as follows. Suppose that $\beta<\alpha$, that $\forceP_\beta$ has been defined, and let $G_\beta$ be a $\forceP_\beta$-generic filter over $W$. For $\gamma<\beta$, write $G_\gamma=G_\beta\upharpoonright\forceP_\gamma$. As in Definition~\ref{0-allowable forcing}, we say that $\eta<\omega_1$ is \emph{free at stage $\beta$} if there are no $\gamma<\beta$ and no finite tuple $\vec z$ such that
    \[
    \dot{\forceQ}_\gamma^{G_\gamma}=\operatorname{Code}(\vec z,\eta).
    \]
    If a proposed coding area $\eta$ is not free, we use instead the least $\eta'<\omega_1$ which is free at stage $\beta$.
    
     \begin{itemize}
    
     \item[(1)] Assume that $F(\beta)=(\dot{m}, \dot{x}, \dot{y},\dot{\eta})$, for a quadruple of $\forceP_{\beta}$-names. We assume that $\dot{x}^{G_{\beta}}=:x$, $\dot{y}^{G_{\beta}}=:y$ are reals, $\dot{m}^{G_{\beta}}=:m$ is a natural number which codes a $\Pi^1_3$-set $A_m$, and $\dot{\eta}^{G_{\beta}}=: \eta$ is an ordinal below $\omega_1$. Let $\eta^*=\eta$ if $\eta$ is free at stage $\beta$, and otherwise let $\eta^*$ be the least coding area which is free at stage $\beta$. We set
     \[
     \dot{\forceQ}_{\beta}^{G_{\beta}}=\forceP(\beta)^{G_{\beta}}:= \operatorname{Code}(m,x,y,\eta^*).
     \]
     Thus this factor writes the exact tuple $(m,x,y)$ into $\vec S$.
    
     \item[(2)] Assume that $F(\beta)=(\dot{m},\dot{x},\dot{z},\dot{\eta})$, where $\dot{x},\dot{z}$ are $\forceP_\beta$-names for reals, $\dot{m}$ is a $\forceP_\beta$-name for the G\"odel number of a $\Sigma^1_5$-formula, and $\dot{\eta}$ is a $\forceP_\beta$-name for an ordinal below $\omega_1$. Let
     \[
     x=\dot{x}^{G_\beta},\qquad z=\dot{z}^{G_\beta},\qquad \eta=\dot{\eta}^{G_\beta}.
     \]
     Let $\eta^*=\eta$ if $\eta$ is free at stage $\beta$, and otherwise let $\eta^*$ be the least coding area which is free at stage $\beta$. We set
     \[
     \dot{\forceQ}_{\beta}^{G_{\beta}}:=\operatorname{Code}(1,1,x,z,\eta^*).
     \]
     Thus this factor writes the exact tuple $(1,1,x,z)$ into $\vec S$. To avoid ambiguity, we assume that $1$ is not the G\"odel number of a formula.

     \end{itemize}

     At every successor stage we set
     \[
     \forceP_{\beta+1}:=\forceP_\beta*\dot{\forceQ}_\beta.
     \]
     If $\lambda<\alpha$ is a nonzero limit ordinal, we take the finite-support direct limit
     \[
     \forceP_\lambda:=\operatorname{dir}\lim_{\beta<\lambda}\forceP_\beta.
     \]
     
        \end{definition}

 This version of ``allowable$"$ is less complicated than the version we needed for the proof of the first theorem. The reason is that we just have two tasks we need to take care of, and that the most tricky case, namely forcing reduction and a failure of uniformization is not existing. 

We proceed now in a very similar way to the proof of the first main theorem, but the thinning out process is tailored to finally obtain a universe where the $\Pi^1_3$-uniformization property is true. We employ the thinning out process as detailed in  \cite{Uniformization} or \cite{Uniformization with wellorder}. Then define a version of ``thinning out while leaving out codes$"$ which is a straightforward adaption of the construction present in this paper.

Equipped with these notions, we start to prove the second main theorem following closely the proof of the first one. As we only have to deal with the $\Pi^1_3$-uniformization and a $\Sigma^1_5$-good wellorder, the proof is shorter.

\subsection{Informal discussion of the construction}\label{InformalSecondMainTheorem}

Before giving the formal thinning-out construction, we describe the two ideas which are combined in the proof.  The first one produces $\Pi^1_3$-uniformizations, while the second one produces a good $\Sigma^1_5$-well-order.  Both use the same coding machinery, but the fixed-point problem occurs only in the first part.

Fix a $\Pi^1_3$-set $A_m\subseteq(\omega^\omega)^2$ and a real $x$.  We want to decide whether the $x$-section of $A_m$ is nonempty and, if it is, to select a single real from it.  The intended definition of the uniformizing function is
\[
f_m(x)=y
\quad\Longleftrightarrow\quad
(m,x,y)\text{ is not coded into }\vec S.
\]
Since the assertion that a tuple is coded into $\vec S$ is $\Sigma^1_3$, this gives a $\Pi^1_3$ definition of the graph.  The forcing construction must therefore arrange that exactly one pair $(m,x,y)$ remains uncoded whenever the $x$-section of $A_m$ is nonempty, and that this distinguished pair belongs to $A_m$ in the final model.

To achieve this, we search for a class $\mathcal P$ of allowable coding forcings for which the following dichotomy is stable.

\begin{enumerate}
\item[(a)] There is a real $y$ such that no forcing in $\mathcal P$ forces
\[
(x,y)\notin A_m.
\]
In this case $(x,y)$ is persistent with respect to the class of forcings used later.  Among the persistent candidates of least rank, we choose the one whose equivalence class of forcing names has the $<$-least canonical code and declare it to be $f_m(x)$.  We then reserve the tuple $(m,x,y)$: the corresponding coding forcing is never used.  Whenever the bookkeeping later presents a real $y'\neq y$, we use a coding forcing for $(m,x,y')$.  Thus $(m,x,y)$ is the unique tuple in the $x$-section which remains uncoded, and its persistence guarantees that $(x,y)$ remains in $A_m$ throughout the later iteration.

\item[(b)] No such persistent real exists.  Equivalently, every candidate $y$ can be forced out of $A_m$ by some forcing in $\mathcal P$.  The bookkeeping presents the possible pairs $(x,y)$ one at a time, and at the corresponding stage we use a forcing in $\mathcal P$ which forces
\[
(x,y)\notin A_m.
\]
Since every relevant name occurs cofinally often, every possible candidate is eventually treated.  In the final model the $x$-section of $A_m$ is therefore empty, and no value of the uniformizing function is required at $x$.
\end{enumerate}

The circular point is that the class $\mathcal P$ used in this dichotomy must itself contain the forcing which is produced when following the above rules. This is the fixed-point problem behind the construction. Starting with the class of all allowable forcings, we thin the class out transfinitely. At a derivative stage, the earlier classes are used only to test whether a persistent candidate has appeared. If none has appeared, the derivative does not insert one of the earlier-rank witness forcings; it merely makes the coding guess prescribed by the bookkeeping. Only after the derivative has stabilized do we use an actual forcing from the fixed-point class to force a bookkeeping pair out of $A_m$. Thus the derivative approximates the above dichotomy, while the fixed-point class makes the forceable-out alternative available without changing the degree of allowability. The ranks recorded in the sets $I_\beta$, together with the keys of the candidates, determine the choice in the persistent case. This is the construction from \cite{Uniformization}; the version combined with the well-order coding is developed in \cite{Uniformization with wellorder}.

The second task is to construct a good $\Sigma^1_5$-well-order of the reals. During the iteration we build an increasing sequence
\[
(<_\beta\mid\beta<\omega_1)
\]
of countable partial well-orders. When the bookkeeping presents reals $x$ and $z$, we extend the current partial order, if necessary, so that $x$ belongs to its field, and use a forcing of the form
\[
\operatorname{Code}(1,1,x,z,\eta)
\]
whenever that exact tuple has not already been forbidden. The bookkeeping returns to every correct pair $(x,z)$ cofinally often. Hence every correct predecessor code is written into $\vec S$ cofinally often.

At the same time, once a real $z'$ is seen not to code the correct predecessor set of some real $x'$ in the partial well-order already constructed, the tuple $(1,1,x',z')$ is added to the forbidden family.  From that stage on, no forcing in the subsequent fixed-point classes is allowed to code this tuple.  Thus a correct predecessor code is coded cofinally often, whereas every incorrect code is eventually excluded and can occur only below some bounded part of the construction.  The relation saying that $z$ codes the set of predecessors of $x$ is consequently $\Sigma^1_5$, and the union
\[
<\;=\bigcup_{\beta<\omega_1}<_\beta
\]
is a good $\Sigma^1_5$-well-order of the reals.

The two parts can be interleaved because a stage devoted to the well-order merely adds new forbidden information and then restarts the thinning-out process inside a subclass of the preceding fixed-point class.  The uniformization decisions already made remain persistent under this restriction.  In the final model, the first part yields $\Pi^1_3$-uniformization and hence $\Sigma^1_4$-uniformization, while the good $\Sigma^1_5$-well-order yields $\Sigma^1_n$-uniformization for every $n\geq5$ by Addison's argument.

We use one rank convention throughout the rest of this section. Each application of the derivative in the main construction belongs to a fixed-point phase indexed by an ordinal $\xi<\omega_1$. If a tentative value becomes persistent at derivative rank $\zeta$ during phase $\xi$, its \emph{global rank} is
\[
\rho=(\xi,\zeta).
\]
Global ranks are ordered lexicographically. Thus earlier fixed-point phases have priority, and within one phase smaller derivative rank has priority. This is compatible with later restarts: a value which is persistent with respect to an earlier fixed-point class remains persistent after that class is restricted. For the first application of the derivative over $W$ we use phase index $\xi=0$. When the phase is fixed, we may still speak of ``$\zeta$-allowable'' forcings for the local derivative class; in the main construction we call this the $(\xi,\zeta)$-allowable class, and we write $(\xi,\infty)$-allowable for its fixed point.

Accordingly, every tentative-value tuple has the form
\[
(m,x,y,\rho,\kappa)
=(m,x,y,(\xi,\zeta),\kappa),
\]
where $\kappa$ is the key assigned when the value is selected.

The product argument requires the candidate search to be localized to a fixed subforcing rather than to the whole initial segment of the iteration. We therefore restore the support notation used in \cite{Uniformization}.

\begin{definition}\label{AdmissibleSupportUniformization}
Let
\[
\forceP=((\forceP_\nu,\dot{\forceQ}_\nu)\mid\nu<\delta)
\]
be an allowable finite-support iteration, and let $\beta\leq\delta$. A set $A\subseteq\beta$ is \emph{admissible for $\forceP_\beta$} if the restricted iteration
\[
\forceP_A:=\Asterisk_{\nu\in A}\dot{\forceQ}^{A}_\nu
\]
taken in increasing order is well defined and is itself allowable, where for every $\nu\in A$ the name $\dot{\forceQ}^{A}_\nu$ is a $\forceP_{A\cap\nu}$-name whose image under the canonical embedding is forced to equal $\dot{\forceQ}_\nu$. We require that the resulting canonical map
\[
\forceP_A\longrightarrow\forceP_\beta
\]
is a complete embedding. If $G_\beta$ is $\forceP_\beta$-generic, we write
\[
G_A:=G_\beta\cap\forceP_A.
\]
An \emph{$A$-local name} is a $\forceP_A$-name.
\end{definition}

We use only admissible supports. Initial segments are the special case $A=\gamma$. Admissible supports are closed under finite unions and under the translations induced by the canonical complete embeddings used below. If $A\subseteq B$ are admissible, then
\[
\forceP_A\lessdot\forceP_B\lessdot\forceP_\beta,
\]
and the corresponding quotients are allowable. These facts follow directly from the recursive definition: every stage belonging to a support depends only on earlier stages belonging to the same support.

Following \cite{Uniformization}, we compare local names for reals through common allowable refinements. If $\forceP,\forceQ$ are allowable forcings, $p\in\forceP$, $q\in\forceQ$, and $\dot a,\dot b$ are names for reals, write
\[
(\forceP,p,\dot a)\sim_0(\forceQ,q,\dot b)
\]
if there are an allowable forcing $\forceR$, complete embeddings
\[
i:\forceP\longrightarrow\forceR,
\qquad
j:\forceQ\longrightarrow\forceR,
\]
whose quotient forcings are allowable, and $r\in\forceR$ such that
\[
r\leq i(p),j(q)
\qquad\text{and}\qquad
r\Vdash_{\forceR}i(\dot a)=j(\dot b).
\]
Let $\sim$ be the equivalence relation generated by $\sim_0$ and define
\[
c(\forceP,p,\dot a)
:=\min_{<}\{\operatorname{code}(\forceQ,q,\dot b):
(\forceQ,q,\dot b)\sim(\forceP,p,\dot a)\}.
\]
Let $N$ be the model over which $\forceP$ is formed and let $G_\beta$ be $\forceP_\beta$-generic over $N$. If $A$ is admissible and $a\in N[G_A]$ is a real, put
\[
\operatorname{key}_A(a)
:=\min_{<}\{c(\forceP_A,p,\dot a):
 p\in G_A,\ \dot a\text{ is a nice $\forceP_A$-name, and }
 \dot a^{G_A}=a\}.
\]
This key does not depend on the admissible support used to present $a$. Indeed, if $A$ and $B$ are admissible and
\[
a\in N[G_A]\cap N[G_B],
\]
then $A\cup B$ is admissible and $\forceP_{A\cup B}$ is a common allowable refinement of $\forceP_A$ and $\forceP_B$. Any two local names which evaluate to $a$ are forced equal below a condition in $G_{A\cup B}$, and hence their presentations are $\sim_0$-equivalent. Thus
\[
\operatorname{key}_A(a)=\operatorname{key}_B(a).
\]
We write $\operatorname{key}(a)$ for this common value. The same argument shows that the key is unchanged when a support and its names are translated along one of the canonical complete embeddings used in a product, an amalgamation, or a later stage block.

Candidates are always compared first by their global rank $\rho$ in the lexicographic order and, among candidates of the same global rank, by their key.

\subsection{The uniformization derivative}\label{UniformizationDerivative}

Fix a fixed-point phase $\xi<\omega_1$ and work in one of the generic extensions $N$ under consideration; for the initial phase we have $N=W$. The definition below is uniform in $N$. At the beginning of the phase we are given a base class of allowable pairs $(\forceP,\dot I)$, where $\dot I$ is the auxiliary tentative-value name. In the initial phase this is the class of the allowable forcings from the preceding definition, equipped with the empty auxiliary name. Later phases start from subclasses obtained by imposing the forbidden information accumulated so far.

For an ordinal $\alpha>0$, assume inductively that the notion of $(\xi,\zeta)$-allowability has been defined for every $\zeta<\alpha$. We define $(\xi,\alpha)$-allowability by modifying the base allowable construction only at uniformization entries of the bookkeeping. The well-order entries are left unchanged.

Let
\[
\forceP=((\forceP_\beta,\dot{\forceQ}_\beta)\mid\beta<\delta)
\]
be a finite-support iteration of length $\delta<\omega_1$. Together with it we recursively construct coherent names
\[
(\dot I_\beta\mid\beta\leq\delta).
\]
Initially
\[
\dot I_0=\check\emptyset,
\qquad
I_0=\emptyset.
\]
Suppose that $\forceP_\beta$ and $\dot I_\beta$ have been defined, let $G_\beta\subseteq\forceP_\beta$ be generic over $N$, and put
\[
I_\beta:=\dot I_\beta^{G_\beta}.
\]
Every local bookkeeping entry used at stage $\beta$ includes an admissible support $A\subseteq\beta$. Its names are $\forceP_A$-names and are interpreted by $G_A$. If a name arises inside a product or an amalgamated allowable block, its support and the name itself are translated together along the corresponding complete embedding.

We distinguish the following cases.

\begin{enumerate}
\item[(a)] \emph{A persistent tentative value appears.}

Assume that the bookkeeping presents a uniformization entry
\[
(A,\dot m,\dot x,\dot y,\dot\eta),
\]
where $A\subseteq\beta$ is admissible for $\forceP_\beta$, the names $\dot m,\dot x,\dot y,\dot\eta$ are $\forceP_A$-names, and, writing
\[
m=\dot m^{G_A},\qquad
x=\dot x^{G_A},\qquad
y=\dot y^{G_A},\qquad
\eta=\dot\eta^{G_A},
\]
we have that $m$ is the G\"odel number of a $\Pi^1_3$-set $A_m$ and $\eta<\omega_1$.

We search for the least $\zeta<\alpha$ for which there is a real $a\in N[G_A]$ such that
\[
N[G_\beta]\models (x,a)\in A_m
\]
and there is no $(\xi,\zeta)$-allowable continuation $\forceR$ of $\forceP_\beta$ such that, in $N[G_\beta]$,
\[
\forceR/G_\beta\Vdash (x,a)\notin A_m.
\]
Here a $(\xi,\zeta)$-allowable continuation means an allowable iteration whose initial segment is $\forceP_\beta$ and whose tail is governed by the $(\xi,\zeta)$-allowable rules computed in the corresponding generic extension. The candidate set is therefore the fixed set of reals $N[G_A]$, while persistence is tested in the full current model $N[G_\beta]$.

If such a $\zeta$ exists, put
\[
\rho:=(\xi,\zeta).
\]
Among all reals $a\in N[G_A]$ witnessing persistence at this least rank, choose $y_0$ with $\operatorname{key}_A(y_0)$ least in the fixed background well-order, and put
\[
\kappa:=\operatorname{key}_A(y_0).
\]
Thus the choice is first by global rank and then by key. The key is stored with the tentative value and is not recomputed in a later, larger universe.

The selected tuple $(m,x,y_0)$ is not coded. If the bookkeeping real $y$ is different from $y_0$, let $\eta^*=\eta$ when $\eta$ is free at stage $\beta$, and otherwise let $\eta^*$ be the least free coding area, and set
\[
\dot{\forceQ}_\beta^{G_\beta}
:=\operatorname{Code}(m,x,y,\eta^*).
\]
If $y=y_0$, choose instead the real $z\in N[G_A]$ with least key among the reals different from $y_0$, let $\dot z$ be the $<$-least nice $\forceP_A$-name for $z$, and use
\[
\dot{\forceQ}_\beta^{G_\beta}
:=\operatorname{Code}(m,x,z,\eta^*)
\]
with $\eta^*$ defined in the same way. In either case the selected tuple $(m,x,y_0)$ is left uncoded.

Let $\dot y_0$ be the $<$-least nice $\forceP_A$-name for $y_0$. After passing all earlier names to $\forceP_{\beta+1}$ by the canonical embeddings, let $\dot I_{\beta+1}$ be the canonical name whose interpretation is
\[
I_{\beta+1}
=I_\beta\cup\{(m,x,y_0,\rho,\kappa)\}
=I_\beta\cup\{(m,x,y_0,(\xi,\zeta),\kappa)\}.
\]
Thus only the tentative value actually selected at this stage is added to the auxiliary set.

\item[(b)] \emph{No persistent tentative value has appeared.}

Assume that the bookkeeping presents a uniformization entry as above and that case~(a) fails. Equivalently, for every $\zeta<\alpha$ and every $a\in N[G_A]$,
\[
N[G_\beta]\models (x,a)\in A_m
\]
implies that there is a $(\xi,\zeta)$-allowable continuation $\forceR$ of $\forceP_\beta$ such that
\[
N[G_\beta]\models
\forceR/G_\beta\Vdash (x,a)\notin A_m.
\]
These witness forcings are used only in deciding that case~(a) fails. They are not inserted into the iteration at this derivative stage.

Instead the derivative makes the bookkeeping guess. Let $\eta^*=\eta$ if $\eta$ is free at stage $\beta$, and otherwise let $\eta^*$ be the least free coding area. We set
\[
\dot{\forceQ}_\beta^{G_\beta}
:=\operatorname{Code}(m,x,y,\eta^*).
\]
No assertion is made here that this coding forcing forces $(x,y)$ out of $A_m$. The role of this clause is only to make the provisional coding decision prescribed by the bookkeeping. No tentative value is selected, so
\[
\dot I_{\beta+1}:=\dot I_\beta,
\qquad
I_{\beta+1}=I_\beta,
\]
after regarding $\dot I_\beta$ as a $\forceP_{\beta+1}$-name.

Only after the derivative hierarchy has reached the fixed point will the main iteration use an actual forcing from that fixed-point class to force a bookkeeping pair $(x,y)$ out of $A_m$.

\item[(c)] \emph{Well-order coding.}

The well-order factors do not participate in the persistence test. If the bookkeeping presents a well-order entry
\[
(A,\dot m,\dot x,\dot z,\dot\eta),
\]
where $A\subseteq\beta$ is admissible, the displayed names are $\forceP_A$-names, $\dot m^{G_A}$ is the G\"odel number of a $\Sigma^1_5$-formula, $\dot x,\dot z$ are names for reals and $\dot\eta$ is a name for a countable ordinal, let
\[
x=\dot x^{G_A},\qquad
z=\dot z^{G_A},\qquad
\eta=\dot\eta^{G_A}.
\]
Let $\eta^*=\eta$ if $\eta$ is free at stage $\beta$, and otherwise let $\eta^*$ be the least free coding area. We use the base allowable rule unchanged and set
\[
\dot{\forceQ}_\beta^{G_\beta}
:=\operatorname{Code}(1,1,x,z,\eta^*).
\]
Thus the derivative codes exactly the pair handed to it by the bookkeeping. Since no tentative uniformizing value is selected here,
\[
\dot I_{\beta+1}:=\dot I_\beta,
\qquad
I_{\beta+1}=I_\beta.
\]
\end{enumerate}

In all other cases we use the trivial forcing and carry $\dot I_\beta$ forward. At every successor stage we put
\[
\forceP_{\beta+1}:=\forceP_\beta*\dot{\forceQ}_\beta.
\]
At a nonzero limit stage $\lambda\leq\delta$ we take the finite-support direct limit
\[
\forceP_\lambda
:=\operatorname{dir}\lim_{\beta<\lambda}\forceP_\beta.
\]
Regard each $\dot I_\beta$, $\beta<\lambda$, as a $\forceP_\lambda$-name via the canonical embeddings and let $\dot I_\lambda$ be the canonical name whose interpretation is their union. Hence
\[
I_\lambda
:=\dot I_\lambda^{G_\lambda}
=\bigcup_{\beta<\lambda}I_\beta.
\]

\begin{definition}\label{UniformizationAlphaAllowable}
Fix a phase $\xi<\omega_1$ and an ordinal $\alpha>0$. A pair $(\forceP,\dot I)$ is called $(\xi,\alpha)$-allowable if it is obtained from the current base class by the preceding recursive rules, where all persistence tests range over the classes $(\xi,\zeta)$ for $\zeta<\alpha$. When $\dot I$ is clear from context we simply call $\forceP$ $(\xi,\alpha)$-allowable. For the initial phase $\xi=0$ we also use the local terminology ``$\alpha$-allowable.''
\end{definition}

The distinction between case~(b) of the derivative and the later forceable-out case is important. The derivative only records which coding decision the bookkeeping guesses when no value has yet become persistent at one of the preceding ranks. Actual force-out forcings are used only after stabilization, when the witness itself belongs to the fixed-point class governing the outer construction.

For a fixed phase $\xi<\omega_1$ and an ordinal $\alpha$, let
\[
\mathcal U_{\xi,\alpha}
\]
denote the class of pairs $(\forceP,\dot I)$ which are $(\xi,\alpha)$-allowable in the sense of Definition~\ref{UniformizationAlphaAllowable}, computed in the generic extension under consideration. At the beginning of a phase, $\mathcal U_{\xi,0}$ denotes the corresponding base class.

\begin{lemma}\label{UniformizationDerivativeProperties}
Fix a phase $\xi<\omega_1$. The classes $(\mathcal U_{\xi,\alpha}\mid\alpha\in\operatorname{Ord})$ have the following properties in every generic extension in which they are computed.
\begin{enumerate}
\item For every $\alpha$, the class $\mathcal U_{\xi,\alpha}$ is nonempty and is uniformly definable from the base class of the phase and the defining parameters of the construction.
\item If $\zeta<\alpha$, then
\[
\mathcal U_{\xi,\alpha}\subseteq\mathcal U_{\xi,\zeta}.
\]
\item Suppose that
\[
(\forceP,\dot I^{\forceP}),\ (\forceQ,\dot I^{\forceQ})
\in\mathcal U_{\xi,\alpha}
\]
and
\[
C^{\forceP}\cap C^{\forceQ}=\emptyset.
\]
Then $\forceP\times\forceQ$ densely embeds into a forcing $\forceR$ for which there is a name $\dot I^{\forceR}$ such that
\[
(\forceR,\dot I^{\forceR})\in\mathcal U_{\xi,\alpha}.
\]
Under the canonical translations of supports and names induced by the product embedding, the auxiliary tentative-value set carried by $\forceR$ is the union of the two auxiliary sets. In particular, every such tuple keeps both the global rank $\rho=(\xi',\zeta')$ and the stored key which it had before the product was formed.
\end{enumerate}
\end{lemma}

\begin{proof}
Nonemptiness and monotonicity are proved by the usual derivative induction. We explain the product clause, which is the reason for using admissible supports.

Proceed simultaneously by induction on $\alpha$ and on the stages of the concatenated iteration. We strengthen the induction statement by requiring that, at every translated second-block entry, both the support-local candidate set and the forceability-out test at each lower derivative rank agree with those in the original second forcing. Let
\[
(\forceP^0,\dot I^0),(\forceP^1,\dot I^1)\in\mathcal U_{\xi,\alpha},
\qquad
\delta_0:=\operatorname{lh}(\forceP^0),
\]
and let $F^0,F^1$ be witnessing bookkeeping functions. Form the concatenated bookkeeping function. In the first block use $F^0$. In the second block, if an entry of $F^1$ at stage $\beta$ carries an admissible support $A\subseteq\beta$, replace it by
\[
A^*:=\{\delta_0+\nu:\nu\in A\}
\]
and translate every $\forceP^1_A$-name to the corresponding $\forceR_{A^*}$-name.

Inductively, before the stage corresponding to $\beta$ the translated support is admissible and
\[
\forceR_{A^*}\cong\forceP^1_A.
\]
Consequently, if $K$ is generic for the concatenated forcing, then
\[
N[K_{A^*}]\cong N[G^1_A].
\]
Thus the reals considered as tentative values in the second block are exactly the reals which were considered in the construction of $\forceP^1$; the generic for $\forceP^0$ does not enlarge this candidate set.

It remains to compare persistence. Fix $a\in N[G^1_A]$ and $\zeta<\alpha$. The strengthened induction hypothesis at rank $\zeta$ identifies the $(\xi,\zeta)$-allowable continuation test over $\forceP^1_\beta$ with the corresponding test over the concatenated prefix. In one direction, a force-out continuation over $\forceP^1_\beta$ is combined with the first block by the lower-rank product property. In the other direction, the support-shifted local presentation supplied by the lower-rank induction is pulled back along
\[
\forceP^1_A\cong\forceR_{A^*}
\]
and the analogous isomorphisms at the later stages of that presentation. Since every local name is carried together with its admissible support, this pullback does not introduce the $\forceP^0$-generic into any local evaluation. The complete embeddings preserve the statement $(x,a)\notin A_m$. Hence $a$ has a $(\xi,\zeta)$-allowable force-out continuation in the original second block if and only if it has one in the concatenated construction.

The least persistence rank is therefore unchanged. The canonical isomorphism
\[
\forceP^1_A\cong\forceR_{A^*}
\]
preserves the key, so the same tentative value is selected and the same coding factor is used. The well-order clause is immediate, since it contains no persistence test. This completes the stage induction and shows that the second block reproduces $\forceP^1$. The auxiliary names are translated and united. Finally, the disjointness of the coding areas allows the concatenation to commute with the first block, giving a dense embedding of $\forceP^0\times\forceP^1$ into the resulting allowable forcing.
\end{proof}

The argument is the support-local form of the uniformization derivative and product induction from \cite{Uniformization}; see also \cite{Uniformization with wellorder} for the simultaneous use of the coding machinery with a definable well-order.

Modulo the fixed coding used throughout, there are only set-many relevant pairs $(\forceP,\dot I)$ of size at most $\aleph_1$. Since the sequence
\[
(\mathcal U_{\xi,\alpha}\mid\alpha\in\operatorname{Ord})
\]
is decreasing and every class is nonempty, it stabilizes. Let $\theta_\xi$ be the least ordinal such that
\[
\mathcal U_{\xi,\theta_\xi}
=
\mathcal U_{\xi,\theta_\xi+1},
\]
and put
\[
\mathcal U_{\xi,\infty}
:=
\mathcal U_{\xi,\theta_\xi}.
\]
We call the members of $\mathcal U_{\xi,\infty}$ the $(\xi,\infty)$-allowable pairs, or, when the auxiliary name is clear, the $(\xi,\infty)$-allowable forcings. By Lemma~\ref{UniformizationDerivativeProperties}, the fixed-point class is nonempty and has the same product property for forcings with disjoint coding areas.

\subsection{Avoidance, products and restarts}
As in the proof of the first main theorem, we use forbidden families of tuples. In this section every forcing name occurring in such a tuple is a local name: it is carried together with an admissible support $A$, although we usually suppress $A$ from the notation. Reals already in the current model are identified with their canonical names on the empty support. If a local name belongs to a tail, a product, or an allowable stage block, we translate both its support and the name itself along the relevant complete embedding. For a finite tuple of local names we take the finite union of their translated supports and regard all components as names over the resulting admissible support.

\begin{definition}\label{AvoidingBUniformization}
Let $B$ be a set of finite tuples of local forcing names satisfying the preceding convention. An allowable iteration
\[
\forceP=((\forceP_\alpha,\dot{\forceQ}_\alpha)\mid\alpha<\delta)
\]
\emph{avoids $B$} if, for every $\alpha<\delta$ and every $\forceP_\alpha$-generic filter $G_\alpha$, there are no $\vec x\in B$, no admissible support $A\subseteq\alpha$ carrying the translated components of $\vec x$, and no $\eta<\omega_1$ such that
\[
\dot{\forceQ}_\alpha^{G_\alpha}
=
\operatorname{Code}(\vec x^{G_A},\eta).
\]
 In particular, if a forbidden tuple has the form
\[
(1,1,x,\dot z),
\]
then the factors excluded by avoidance are exactly the forcings
\[
\operatorname{Code}(1,1,x,z,\eta)
\]
for the corresponding interpretation $z$ and coding area $\eta$.
\end{definition}

We next carry out the derivative with a fixed forbidden family. This is the version used at the beginning of the construction and, after a restart, with the current fixed-point class in place of the original base class.

\begin{definition}\label{FixedBUniformizationDerivative}
Fix a phase $\xi<\omega_1$ and a forbidden family $B$ as above. Let
\[
\mathcal U^B_{\xi,0}
:=
\{(\forceP,\dot I):\forceP\text{ is $0$-allowable, }\dot I=\check\emptyset,
\text{ and }\forceP\text{ avoids }B\}.
\]
For $\alpha>0$, define $\mathcal U^B_{\xi,\alpha}$ by the uniformization derivative from Definition~\ref{UniformizationAlphaAllowable}, starting from $\mathcal U^B_{\xi,0}$. Thus every persistence test at derivative rank $\alpha$ is made against the classes
\[
\mathcal U^B_{\xi,\zeta}
\qquad(\zeta<\alpha),
\]
every actual coding factor is required to avoid $B$, and every tentative value selected for the first time at derivative rank $\zeta$ is assigned the global rank $(\xi,\zeta)$. The auxiliary names, their admissible supports, and their translations are updated exactly as in the derivative above.
\end{definition}

The product property remains valid for these classes.

\begin{lemma}\label{ProductLemmaUniformizationAvoidingB}
Let $B$ be fixed, let $\alpha$ be an ordinal, and suppose that
\[
(\mathbb P^0,\dot I^0),(\mathbb P^1,\dot I^1)
\in\mathcal U^B_{\xi,\alpha}
\]
with
\[
C^{\mathbb P^0}\cap C^{\mathbb P^1}=\emptyset.
\]
Then $\mathbb P^0\times\mathbb P^1$ densely embeds into a forcing $\mathbb R$ for which there is a name $\dot I^{\mathbb R}$ such that
\[
(\mathbb R,\dot I^{\mathbb R})
\in\mathcal U^B_{\xi,\alpha}.
\]
Under the canonical translations of supports and names induced by the product embedding, the interpretation of $\dot I^{\mathbb R}$ is the union of the interpretations of $\dot I^0$ and $\dot I^1$. In particular, every tentative value keeps its previously assigned global rank and stored key. The same statement holds at the fixed point.
\end{lemma}

\begin{proof}
Repeat the proof of Lemma~\ref{UniformizationDerivativeProperties}. In the second block, every admissible support $A$ is replaced by its shifted copy $A^*$, so the corresponding support extension and candidate set are unchanged. The lower-rank product induction gives the same persistence decisions, and hence every actual factor in the amalgamated forcing is a translated factor from one of the two original forcings. The translated local names still denote the same tuples. Since both original forcings avoid $B$, the amalgamated forcing avoids $B$ as well. The fixed-point statement follows at the stabilizing rank.
\end{proof}

As before, the disjointness hypothesis causes no restriction in applications, since the countably many coding areas used by one of the forcings may be moved to fresh coordinates before the product is formed.

The classes
\[
(\mathcal U^B_{\xi,\alpha}\mid\alpha\in\operatorname{Ord})
\]
are nonempty and decreasing. By the same set-sized stabilization argument used above, there is a least $\theta^B_\xi$ such that
\[
\mathcal U^B_{\xi,\theta^B_\xi}
=
\mathcal U^B_{\xi,\theta^B_\xi+1}.
\]
We put
\[
\mathcal U^B_{\xi,\infty}
:=
\mathcal U^B_{\xi,\theta^B_\xi}
\]
and call its members the $(\xi,\infty)$-allowable pairs avoiding $B$.

\subsubsection{Restarting the derivative}

During the main construction the forbidden family is enlarged only after a stage forcing has been used. We therefore need a restart operation which begins inside the current fixed-point class rather than again from the class of all $0$-allowable forcings.

\begin{definition}\label{RestartUniformizationDerivative}
Work in one of the generic extensions under consideration and let $\mathcal C_\infty$ be the current nonempty fixed-point class of allowable pairs. Assume that $\mathcal C_\infty$ has the support-preserving product property for forcings with disjoint coding areas. Let $D$ be the newly added forbidden family, satisfying the exact-tuple and local-support convention from Definition~\ref{AvoidingBUniformization}. Earlier forbidden information is already incorporated into $\mathcal C_\infty$.

Define
\[
\operatorname{Rest}(\mathcal C_\infty,D)_0
:=
\{(\forceP,\dot I)\in\mathcal C_\infty:
\forceP\text{ avoids }D\}.
\]
Fix the phase index $\xi$ assigned to this restart. For $\alpha>0$, define $\operatorname{Rest}(\mathcal C_\infty,D)_\alpha$ by applying the uniformization derivative to the base class $\operatorname{Rest}(\mathcal C_\infty,D)_0$. Thus persistence at derivative rank $\alpha$ is tested only against the preceding classes
\[
\operatorname{Rest}(\mathcal C_\infty,D)_\zeta
\qquad(\zeta<\alpha),
\]
and every new tentative value first selected at derivative rank $\zeta$ receives the global rank $(\xi,\zeta)$. The auxiliary tentative-value information already carried by a member of $\mathcal C_\infty$ is retained, and newly selected tentative values are adjoined to it.

The resulting sequence is nonempty and decreasing, has the support-preserving product property on disjoint coding areas, and stabilizes by the same argument as above. We write
\[
\operatorname{Rest}(\mathcal C_\infty,D)_\infty
\]
for its fixed point. We call the passage
\[
\mathcal C_\infty\longmapsto\operatorname{Rest}(\mathcal C_\infty,D)_\infty
\]
the \emph{restart of the uniformization derivative with the new forbidden family $D$}.
\end{definition}

The product assertion in the definition is proved by the same induction as Lemma~\ref{ProductLemmaUniformizationAvoidingB}: in the second block one shifts each admissible support together with its names. Hence the candidate set remains the same, and the exact-tuple avoidance requirement is preserved.

Since
\[
\operatorname{Rest}(\mathcal C_\infty,D)_0\subseteq\mathcal C_\infty,
\]
every class produced by the restart is a subclass of the class which governed the construction before the restart. Thus all restrictions imposed at earlier phases remain in force; the definition of the restart needs to mention only the newly added family $D$.

If no new forbidden information is added at a stage, no restart is performed. In that case the construction continues with the current fixed-point class, reinterpreted in the new generic extension. The forcing stages and the successive restarts will be defined simultaneously when the main iteration is specified.

\subsection{The closure operation}

The forceable-out stages of the main iteration use a closure of a forcing from the current fixed-point class similar to the closure notion in the proof of the first main theorem. As there this operation is applied only after the derivative has reached its fixed point.

\begin{definition}\label{UniformizationClosureDefinition}
Work in a current model and let $\mathcal C_\infty$ be the current fixed-point class, with its current forbidden family already incorporated. Let
\[
(\forceP,\dot I)\in\mathcal C_\infty,
\qquad
\forceP=((\forceP_\nu,\dot{\forceQ}_\nu)\mid\nu<\delta).
\]
Define $C^1(\forceP)$ by inspecting every actual uniformization coding factor
\[
\dot{\forceQ}_\nu^{G_\nu}
=\operatorname{Code}(m,x,y,\eta).
\]
All local names which produced this factor retain their admissible supports. In the corresponding intermediate model, if some member of $\mathcal C_\infty$ forces
\[
(x,y)\notin A_m,
\]
choose the least such witness, move its coding areas away from those of the current factor, and use the support-preserving product property to amalgamate it with that factor. The complete embedding translates every local support and name together. If no such witness exists, leave the factor unchanged. The auxiliary tentative-value name for $C^1(\forceP)$ is the translated union of $\dot I$ with the auxiliary names carried by the inserted witness forcings; all previously assigned global ranks and stored keys are kept unchanged.

Set
\[
C^0(\forceP):=\forceP,
\qquad
C^{n+1}(\forceP):=C^1(C^n(\forceP)),
\]
and let
\[
C^\omega(\forceP)
:=\operatorname{dir}\lim_{n<\omega}C^n(\forceP)
\]
with finite support. Its auxiliary tentative-value name is the translated union of the auxiliary names occurring at the finite stages.
\end{definition}

\begin{lemma}\label{UniformizationClosureLemma}
Let $(\forceP,\dot I)\in\mathcal C_\infty$. Then:
\begin{enumerate}
\item $C^\omega(\forceP)$ belongs to the same current fixed-point class $\mathcal C_\infty$;
\item it avoids the same forbidden family;
\item its auxiliary tentative-value set is the union of the tentative-value tuples carried by $\forceP$ and by all witness forcings inserted during the closure, with all global ranks and stored keys unchanged;
\item if an actual factor of $C^\omega(\forceP)$ codes $(m,x,y)$ and, when that factor appears, some forcing in the current fixed-point class forces $(x,y)\notin A_m$, then
\[
C^\omega(\forceP)\Vdash(x,y)\notin A_m.
\]
\end{enumerate}
\end{lemma}

\begin{proof}
At each finite step, the fresh-coding-area argument and the support-preserving product property keep the replacement factor in $\mathcal C_\infty$ and preserve avoidance. Every admissible support and local name is carried through the corresponding complete embedding. Hence every $C^n(\forceP)$ belongs to the same class, and so does the finite-support direct limit. Finally, every coding factor occurring at some finite stage is treated at the next closure step whenever a force-out witness exists. This proves~(4).
\end{proof}

\begin{lemma}\label{UniformizationPersistenceTailLemma}
Suppose that the tuple
\[
(m,x,y,\rho,\kappa),\qquad \rho=(\xi,\zeta),
\]
is inserted into the auxiliary tentative-value set because
\[
(x,y)\in A_m
\]
and no forcing in the class against which persistence is tested at rank $\rho$ forces
\[
(x,y)\notin A_m.
\]
Then the following hold throughout the later main construction.
\begin{enumerate}
\item Every later fixed-point class is a subclass of the class against which the persistence of $y$ was tested.
\item Every later stage block and every bounded quotient of the main iteration belongs, after translating each admissible support together with its local names, to the corresponding earlier class. This remains true for products and amalgamations, for the closures $C^\omega(\forceP)$, and after later restarts.
\item No later bounded quotient forces $(x,y)\notin A_m$.
\item Consequently,
\[
(x,y)\in A_m
\]
in the final model.
\end{enumerate}
\end{lemma}

\begin{proof}
This is the tail-persistence argument from \cite{Uniformization}; see also \cite{Uniformization with wellorder}. Monotonicity of the derivative and Definition~\ref{RestartUniformizationDerivative} give~(1). For~(2), subsequent stage blocks are chosen from the current fixed-point class and hence from the earlier class. Products stay in that class by the support-preserving Lemma~\ref{ProductLemmaUniformizationAvoidingB}, and closures stay in it by Lemma~\ref{UniformizationClosureLemma}. A restart only passes to a subclass. The admissible supports attached to local names are translated along all these embeddings, so the candidate universe used in the original persistence decision is never enlarged. Thus a bounded quotient forcing $(x,y)\notin A_m$ would contradict the persistence test which introduced $(m,x,y,\rho,\kappa)$, proving~(3). The standard $\Sigma^1_3$ persistence argument then gives~(4).
\end{proof}

\subsection{The main iteration}\label{SecondMainIteration}

We now define the main $\omega_1$-length iteration over $W$. The forcing is constructed simultaneously with the auxiliary tentative-value information, the partial well-order, the new forbidden information, and the current derivative classes. For $\beta<\omega_1$ we recursively define
\[
\forceP_\beta,\qquad G_\beta,\qquad M_\beta=W[G_\beta],\qquad
\dot I_\beta,\qquad I_\beta=\dot I_\beta^{G_\beta},\qquad <_\beta,
\]
together with the family $B_\beta$ of new forbidden tuples produced at stage $\beta$ and the current classes
\[
\mathcal U_{\beta,\alpha}\qquad(\alpha\in\operatorname{Ord}),
\qquad
\mathcal U_{\beta,\infty}.
\]
Here $B_\beta$ records only the new forbidden information produced at stage $\beta$; all earlier restrictions are already incorporated into the current fixed-point class.

\medskip
\noindent
\emph{Initial stage.}
Set
\[
\forceP_0=\mathbf 1,\qquad
G_0=\emptyset,\qquad
M_0=W,\qquad
\dot I_0=\check\emptyset,\qquad
I_0=\emptyset,\qquad
<_0=\emptyset.
\]
No forbidden information has yet been imposed. Let
\[
\mathcal U_{0,0}
:=
\{(\forceR,\check\emptyset):
\forceR\text{ is $0$-allowable for the present construction}\}.
\]
Apply the derivative from Subsection~\ref{UniformizationDerivative} with phase index $0$, and let
\[
\mathcal U_{0,\infty}
\]
be the resulting fixed-point class.

\medskip
\noindent
\emph{Bookkeeping.}
Fix a tagged bookkeeping function
\[
F:\omega_1\longrightarrow H(\omega_1)
\]
such that every relevant tagged tuple of outer names occurs unboundedly often. A uniformization entry has the form
\[
\langle\mathrm{unif},\gamma,\dot m,\dot x,\dot y\rangle,
\]
and a well-order entry has the form
\[
\langle\mathrm{wo},\gamma,\dot x,\dot z\rangle,
\]
where the displayed names are $\forceP_\gamma$-names and $\gamma\leq\beta$ when the entry is used at stage $\beta$. We interpret them with
\[
G_\gamma=G_\beta\cap\forceP_\gamma.
\]
The ordinal $\gamma$ only records an outer stage at which the names are available. It does not replace the admissible supports used inside an allowable stage block. If a name or tentative-value tuple arises inside such a block, a product, or an amalgamation, its admissible support, its local name, its global rank and its stored key are retained and translated along the relevant complete embedding. The bookkeeping returns unboundedly often to every outer name obtained in this way.

\medskip
\noindent
\emph{Successor stage.}
Suppose that the construction has been carried out through stage $\beta$ and that
\[
M_\beta=W[G_\beta]
\]
contains the current fixed-point class $\mathcal U_{\beta,\infty}$. According to the tag $F(\beta)$, one of the cases below supplies a $\forceP_\beta$-name $\dot{\forceQ}_\beta$ and the corresponding auxiliary tentative-value information so that, suppressing the auxiliary name from the notation,
\[
\dot{\forceQ}_\beta^{G_\beta}\in\mathcal U_{\beta,\infty}.
\]
Define
\[
\forceP_{\beta+1}:=\forceP_\beta*\dot{\forceQ}_\beta.
\]
If $H_\beta$ is $\dot{\forceQ}_\beta^{G_\beta}$-generic over $M_\beta$, put
\[
G_{\beta+1}:=G_\beta*H_\beta,
\qquad
M_{\beta+1}:=M_\beta[H_\beta].
\]
The same case specifies the $\forceP_{\beta+1}$-name $\dot I_{\beta+1}$, the next partial well-order $<_{\beta+1}$, and the new forbidden family $B_\beta$.

If $B_\beta=\emptyset$, no restart is performed. We reinterpret the current classes in $M_{\beta+1}$ and put
\[
\mathcal U_{\beta+1,\alpha}
:=(\mathcal U_{\beta,\alpha})^{M_{\beta+1}}
\qquad(\alpha\in\operatorname{Ord}),
\]
and
\[
\mathcal U_{\beta+1,\infty}
:=(\mathcal U_{\beta,\infty})^{M_{\beta+1}}.
\]
If $B_\beta\neq\emptyset$, restart the derivative in $M_{\beta+1}$ inside the current fixed-point class, with phase index $\beta+1$, and put
\[
\mathcal U_{\beta+1,\alpha}
:=
\operatorname{Rest}\bigl((\mathcal U_{\beta,\infty})^{M_{\beta+1}},B_\beta\bigr)_\alpha
\]
for every $\alpha$, together with the corresponding fixed point
\[
\mathcal U_{\beta+1,\infty}
:=
\operatorname{Rest}\bigl((\mathcal U_{\beta,\infty})^{M_{\beta+1}},B_\beta\bigr)_\infty.
\]
All earlier restrictions remain built into this class.

A stage forcing may itself be a countable finite-support allowable iteration. We use it as one block of the outer iteration; its internal coding factors, admissible supports and auxiliary names remain internal to that block.

\medskip
\noindent
\emph{Uniformization stage: a persistent value is available.}
Suppose
\[
F(\beta)=\langle\mathrm{unif},\gamma,\dot m,\dot x,\dot y\rangle,
\qquad \gamma\leq\beta,
\]
and put
\[
m=\dot m^{G_\gamma},\qquad
x=\dot x^{G_\gamma},\qquad
y=\dot y^{G_\gamma}.
\]
Assume that $m$ is the G\"odel number of the $\Pi^1_3$-set $A_m$. For a real $a\in M_\beta$, a global rank $\rho$, and a key $\kappa$, say that $a$ is a \emph{current tentative value} for $(m,x)$, with rank $\rho$ and key $\kappa$, if
\[
(x,a)\in A_m
\]
and there are $(\forceR,\dot J)\in\mathcal U_{\beta,\infty}$ and $p\in\forceR$ such that
\[
p\Vdash_{\forceR}(m,\check x,\check a,\check\rho,\check\kappa)\in\dot J.
\]
If such tentative values exist, let $(\rho_0,\kappa_0)$ be the least pair in the lexicographic order, and let $y_0$ be the corresponding real; if necessary, use the least supported presentation witnessing this pair. By the way these tentative values are introduced by the derivative and by Lemma~\ref{UniformizationPersistenceTailLemma}, $y_0$ is persistent with respect to the current fixed-point class.

If $y\neq y_0$, let $\eta$ be the least coding area which is free at the present stage and put
\[
\dot{\forceQ}_\beta^{G_\beta}
:=\operatorname{Code}(m,x,y,\eta).
\]
If $y=y_0$, use the trivial forcing. Thus the selected tuple $(m,x,y_0)$ is left uncoded at this stage. Regard the earlier names as $\forceP_{\beta+1}$-names through the canonical embeddings and let $\dot I_{\beta+1}$ be the canonical name whose interpretation is
\[
I_{\beta+1}
=I_\beta\cup\{(m,x,y_0,\rho_0,\kappa_0)\}.
\]
Finally set
\[
<_{\beta+1}:=<_\beta,
\qquad
B_\beta:=\emptyset.
\]

\medskip
\noindent
\emph{Uniformization stage: no persistent value is available.}
Suppose the same tagged entry is presented, but the current fixed-point class carries no tentative value for $(m,x)$. Since the derivative has stabilized, the fixed-point dichotomy implies that the bookkeeping candidate $y$ can be forced out of $A_m$ by a member of the current fixed-point class. Choose the least pair
\[
(\forceR_{m,x,y},\dot J_{m,x,y})\in\mathcal U_{\beta,\infty}
\]
such that
\[
\forceR_{m,x,y}\Vdash(x,y)\notin A_m.
\]
If $(x,y)\notin A_m$ already holds in $M_\beta$, the trivial forcing may serve as this witness. Use the closure of the witnessing forcing as the stage forcing:
\[
\dot{\forceQ}_\beta^{G_\beta}
:=C^\omega(\forceR_{m,x,y}).
\]
Let $\dot J_\beta$ be the auxiliary tentative-value name carried by this closure, with all its local supports, ranks and keys translated to the stage forcing. Let $\dot I_{\beta+1}$ be the canonical name whose interpretation is
\[
I_{\beta+1}=I_\beta\cup J_\beta,
\qquad
J_\beta:=\dot J_\beta^{G_{\beta+1}}.
\]
Set
\[
<_{\beta+1}:=<_\beta,
\qquad
B_\beta:=\emptyset.
\]

\medskip
\noindent
\emph{Well-order stage.}
Suppose
\[
F(\beta)=\langle\mathrm{wo},\gamma,\dot x,\dot z\rangle,
\qquad \gamma\leq\beta,
\]
and put
\[
x=\dot x^{G_\gamma},\qquad
z=\dot z^{G_\gamma}.
\]
If $x\in\operatorname{field}(<_\beta)$, let
\[
<_{\beta+1}:=<_\beta.
\]
Otherwise add $x$ as a new greatest element:
\[
<_{\beta+1}
:=<_\beta\cup\{(u,x):u\in\operatorname{field}(<_\beta)\}.
\]
Thus every successor step is an end-extension.

If the exact tuple $(1,1,x,z)$ is not already excluded by the current fixed-point class, let $\eta$ be the least free coding area and put
\[
\dot{\forceQ}_\beta^{G_\beta}
:=\operatorname{Code}(1,1,x,z,\eta).
\]
If that tuple is already forbidden, use the trivial forcing. No tentative value is selected at a pure well-order stage, so
\[
\dot I_{\beta+1}:=\dot I_\beta,
\qquad
I_{\beta+1}=I_\beta.
\]

In $M_{\beta+1}$ let $B_\beta$ consist of the newly forbidden exact tuples
\[
(1,1,x',\dot z')
\]
where $x'\in\operatorname{field}(<_{\beta+1})$, the name $\dot z'$ is carried by an admissible support in an allowable continuation from the current class, and that continuation forces
\[
\{(\dot z')_n:n\in\omega\}
\neq
\{u:u<_{\beta+1}x'\}.
\]
Tuples already excluded by the current fixed-point class are omitted from $B_\beta$. We then restart the derivative with this new family. An incorrect bookkeeping tuple may be coded at the present stage, but after the restart it cannot be coded again. Correct predecessor codes are never put into a forbidden family.

\medskip
\noindent
\emph{Limit stage.}
Let $0<\lambda<\omega_1$ be a limit ordinal and suppose that the construction has been carried out below $\lambda$. Take the finite-support direct limit
\[
\forceP_\lambda
:=\operatorname{dir}\lim_{\beta<\lambda}\forceP_\beta,
\qquad
G_\lambda:=\bigcup_{\beta<\lambda}G_\beta,
\qquad
M_\lambda:=W[G_\lambda].
\]
Regard each $\dot I_\beta$, $\beta<\lambda$, as a $\forceP_\lambda$-name and let $\dot I_\lambda$ be the canonical name whose interpretation is
\[
I_\lambda=\bigcup_{\beta<\lambda}I_\beta.
\]
Since the partial well-orders are extended by end-extension, put
\[
<_\lambda:=\bigcup_{\beta<\lambda}<_\beta.
\]

In $M_\lambda$ reinterpret all preceding fixed-point classes and set
\[
\mathcal U_{\lambda,0}
:=
\bigcap_{\beta<\lambda}(\mathcal U_{\beta,\infty})^{M_\lambda}.
\]
The trivial pair belongs to this intersection. If two members use disjoint coding areas, concatenate their bookkeeping functions and shift every admissible support in the second block together with its names, as in Lemma~\ref{UniformizationDerivativeProperties}. The same local candidate sets, persistence decisions, ranks and stored keys are obtained in every preceding phase. Hence the resulting pair belongs to every preceding fixed-point class and therefore to the intersection. Applying the derivative with phase index $\lambda$ yields a nonempty fixed point
\[
\mathcal U_{\lambda,\infty}
\]
with the support-preserving product property.

Finally define
\[
\forceP_{\omega_1}
:=\operatorname{dir}\lim_{\beta<\omega_1}\forceP_\beta
\]
with finite support. If $G_{\omega_1}\subseteq\forceP_{\omega_1}$ is generic over $W$, then
\[
G_{\omega_1}=\bigcup_{\beta<\omega_1}G_\beta,
\qquad
I_{\omega_1}=\bigcup_{\beta<\omega_1}I_\beta,
\qquad
<_{\omega_1}=\bigcup_{\beta<\omega_1}<_\beta.
\]
Every stage forcing is ccc and has size at most $\aleph_1$. Hence the finite-support iteration $\forceP_{\omega_1}$ is ccc. Since $W\models\CH$ and the iteration has length $\omega_1$, it has size at most $\aleph_1$. Thus
\[
W[G_{\omega_1}]\models\CH.
\]

\begin{lemma}\label{NoUnwantedCodesSecondMainTheorem}
The no-unwanted-codes lemma applies to the allowable classes used in this section. In particular, in the final model a finite tuple $\vec z$ is coded into $\vec S$ if and only if an actual factor
\[
\operatorname{Code}(\vec z,\eta)
\]
occurs in the main iteration or in one of its allowable stage blocks.
\end{lemma}

\begin{proof}
All forcings used here are iterations of the coding forcings from Section~3. The support subforcings, auxiliary sets, closure operation and forbidden families introduce no new kind of coding factor. Hence Lemma~\ref{nounwantedcodes} applies at every bounded stage. Any real witnessing a code in the full $\omega_1$-iteration appears at a bounded stage, so the same conclusion follows.
\end{proof}

\subsection{Verification of the second main theorem}

Let $\forceP_{\omega_1}$ denote the forcing obtained by carrying the preceding recursion through $\omega_1$ stages with the tagged bookkeeping function $F$.

\begin{theorem}
Let $G_{\omega_1}$ be a $\forceP_{\omega_1}$-generic filter over $W$. Then in $W[G_{\omega_1}]$ the $\Pi^1_3$-uniformization property holds and $\Sigma^1_n$-uniformization holds for every $n\geq4$. Moreover, the reals admit a good $\Sigma^1_5$-well-order.
\end{theorem}

\begin{proof}
For a tuple $(m,x,y)$ let
\[
\operatorname{Coded}_{\vec S}(m,x,y)
\]
abbreviate the $\Sigma^1_3$-statement from Section~3 saying that $(m,x,y)$ is coded into $\vec S$. If $m$ is the G\"odel number of a $\Pi^1_3$-set $A_m$ in the plane, define
\[
f_m(x)=y
\quad\Longleftrightarrow\quad
(x,y)\in A_m
\ \text{ and }\
\neg\operatorname{Coded}_{\vec S}(m,x,y).
\]
This is a $\Pi^1_3$ definition.

Fix $m$ and a real $x$ in the final model. We distinguish two cases.

\medskip
\noindent
\emph{A tentative value occurs.}
Assume that
\[
(m,x,y,\rho,\kappa)\in I_{\omega_1}.
\]
Choose the least pair $(\rho_0,\kappa_0)$ occurring for $(m,x)$ and let $y_0$ be the corresponding real, using the least supported presentation in case of a tie. We claim first that $(m,x,y_0)$ was not coded before $y_0$ became the selected value.

Suppose otherwise and consider the first actual factor which coded $(m,x,y_0)$. There are two possibilities. If it was produced directly by a persistent-value stage of the outer iteration, then the value selected at that stage had a smaller rank--key pair than $(\rho_0,\kappa_0)$; if the selected value was $y_0$ itself, the stage did not code $(m,x,y_0)$. The first alternative contradicts the minimality of $(\rho_0,\kappa_0)$ and the second contradicts the choice of the factor.

It remains to consider a factor occurring inside a closed forceable-out block. Let $N_0$ be the model over which that allowable block is computed and let $A$ be the admissible support of the local bookkeeping entry which produced the factor. Then $y_0\in N_0[G_A]$. If a forcing in the current fixed-point class could force $(x,y_0)\notin A_m$, Lemma~\ref{UniformizationClosureLemma} would make the closed block force
\[
(x,y_0)\notin A_m,
\]
contrary to Lemma~\ref{UniformizationPersistenceTailLemma}. If no such force-out forcing exists, then, at the fixed point, one further application of the support-local derivative selects a persistent value from the same candidate universe $N_0[G_A]$. It either selects $y_0$, in which case $(m,x,y_0)$ is left uncoded, or selects a tentative value with a smaller rank--key pair, which is carried into $I_{\omega_1}$. Both alternatives are impossible. Thus no factor codes $(m,x,y_0)$ before it is selected.

After $y_0$ has been selected, Lemma~\ref{UniformizationPersistenceTailLemma} gives
\[
(x,y_0)\in A_m
\]
in the final model. Every later fixed-point class is a subclass of the class against which its persistence was tested. Products, closures and restarts translate the admissible support of its local presentation and preserve the stored pair $(\rho_0,\kappa_0)$; they do not repeat the candidate search in a larger support extension. Since this is the least pair which occurs in $I_{\omega_1}$ for $(m,x)$, the rank/key convention continues to select $y_0$ whenever $(m,x)$ is considered. Hence no later factor codes $(m,x,y_0)$. By Lemma~\ref{NoUnwantedCodesSecondMainTheorem},
\[
\neg\operatorname{Coded}_{\vec S}(m,x,y_0).
\]
Thus $f_m(x)=y_0$.

Let $y'\neq y_0$ be any real in the final model. After names for $m,x,y'$ have appeared, the tagged bookkeeping presents them at a later uniformization stage. The tentative value with least rank--key pair remains $y_0$, with rank $\rho_0$ and key $\kappa_0$, so the stage uses
\[
\operatorname{Code}(m,x,y',\eta)
\]
for a free coding area. Hence $y_0$ is the unique real satisfying the definition of $f_m(x)$.

\medskip
\noindent
\emph{No tentative value occurs.}
Assume that $I_{\omega_1}$ contains no tuple of the form
\[
(m,x,y,\rho,\kappa).
\]
Let $y$ be any real in the final model and choose a later uniformization stage after names for $m,x,y$ have appeared. If $(x,y)\notin A_m$ already holds at that stage, it remains true in the final model, since this is a $\Sigma^1_3$-statement. Otherwise the current fixed-point class carries no persistent tentative value for $(m,x)$. The forceable-out clause therefore chooses a forcing $\forceR_{m,x,y}$ with
\[
\forceR_{m,x,y}\Vdash(x,y)\notin A_m
\]
and uses
\[
C^\omega(\forceR_{m,x,y})
\]
as the stage forcing. Hence $(x,y)\notin A_m$ holds from that stage on. Since every real $y$ is eventually treated,
\[
(A_m)_x=\emptyset
\]
in the final model.

We have proved that $f_m$ has exactly one value on every nonempty section of $A_m$ and no value on an empty section. Since $m$ was arbitrary, $\Pi^1_3$-uniformization holds.

We next verify the well-order conclusion. Let
\[
<\;:=\bigcup_{\beta<\omega_1}<_\beta.
\]
Every successor step is an end-extension and every proper initial segment is countable. Every real in the final model is the interpretation of a name from some bounded outer stage, and the tagged bookkeeping presents that name later at a well-order stage. Hence every real enters the field of $<$. Since the final model satisfies $\CH$, the order has type $\omega_1$.

If $z$ codes the set of $<$-predecessors of $x$, then after $x$ enters the order this predecessor set never changes. The exact tuple $(1,1,x,z)$ is never put into a forbidden family, and the bookkeeping codes it cofinally often. Conversely, if $z$ does not code the correct predecessor set, choose a later well-order stage at which a supported name for $z$ is considered after $x$ has entered the order. The tuple $(1,1,x,z)$ is then added to the new forbidden family. It can be coded only below a bounded stage and is therefore not coded cofinally often. Consequently
\[
z\text{ codes the set of $<$-predecessors of }x
\quad\Longleftrightarrow\quad
(1,1,x,z)\text{ is coded cofinally often}.
\]
The right-hand side is a $\Sigma^1_5$-formula, as in the proof of the first main theorem. Thus $<$ is a good $\Sigma^1_5$-well-order.

Finally, $\Pi^1_3$-uniformization implies $\Sigma^1_4$-uniformization by uniformizing a $\Pi^1_3$ relation on pairs and projecting the selected pair. The good $\Sigma^1_5$-well-order gives $\Sigma^1_n$-uniformization for every $n\geq5$ by Addison's argument. Hence
\[
\Sigma^1_n\text{-uniformization holds for every }n\geq4.
\]
\end{proof}

    \section{Open questions}

    We end with several questions which are related to this article.
    The first observation is concerned with a hypothetical forcing $\forceP$ which would force the $\Pi^1_5$-uniformization property over $L$.
    The construction of such a forcing $\forceP$ necessarily cannot be combined with the methods presented here, as $\Sigma^1_5$ and $\Pi^1_5$-uniformization contradict each other. Given the flexibility of the way the good $\Sigma^1_5$-wellorder is forced here, which can even be arranged such that one can force towards $\Pi^1_3$-uniformization, it is not clear at all how such a hypothetical forcing would look like. The present work seems to suggest that, if it can be done at all, a very new approach is in need.
    \begin{question}
   Can one force $\Pi^1_5$-uniformization over just $L$?
    \end{question}
    A second question concerns patterns of the $\Sigma$-uniformization property. 
        \begin{question}
    Given a real $r \in 2^{\omega}$, can one force a universe where $\Sigma^1_n$-uniformization is true whenever $r(n)=1$ and $\Sigma^1_n$-uniformization fails whenever $r(n)=0$?
    \end{question}
  
    By a classical result of Novikov, for any projective pointclass $\Gamma$, it is impossible to have $\Gamma$ and $\check{\Gamma}$-reduction simultaneously. The case  for separation is still unknown.
    
    \begin{question}
    
    Can one force a universe where $\Sigma^1_3$- and $\Pi^1_3$-separation hold simultaneously?
    
    \end{question}
    
    \section{Acknowledgment}
 
 To the one whose very short life prompted this work. Your family loves you.
 This research was funded in whole by the Austrian Science Fund (FWF) Grant-DOI 10.55776/P37228.


\begin{thebibliography}{20}
    
    
    \bibitem{Addison} J. Addison \textit{Some consequences of the axiom of constructibility}, Fundamenta Mathematica, vol. 46 (1959), pp. 337–357.
    
    
    
    \bibitem{David}
    
    R. David \textit{A very absolute $\Pi^1_2$-real singleton}. Annals of Mathematical Logic 23, pp. 101-120, 1982.
    
    \bibitem{NS saturated and definable} S. Hoffelner \textit{$\hbox{NS}_{\omega_1}$ $\Delta_1$-definable and saturated.}  Journal of Symbolic Logic 86 (1), pp. 25 - 59, 2021.
   
  
    \bibitem{Separation} S. Hoffelner \textit{Forcing the $\Sigma^1_3$-separation property}. Journal of Mathematical Logic 22, No. 2, 2022.
    

     \bibitem{Reduction} S. Hoffelner \textit{Forcing the $\Pi^1_3$-Reduction Property and a Failure of $\Pi^1_3$-Uniformization},  Annals of Pure and Applied Logic 174 (8), 2023.
      
    \bibitem{Forcing axioms and uniformization} S. Hoffelner \textit{Forcing Axioms and the Uniformization Property} Annals of Pure and Applied Logic 175 (10), 2024.
    
    
     \bibitem{BPFA and uniformization} S. Hoffelner \textit{The global $\Sigma^1_{n+2}$-Uniformization Property and $\BPFA$} Advances in Mathematics 470, 2025.

     \bibitem{MA and failure of separation} S. Hoffelner \textit{$\mathsf{MA} (\mathcal{I})$ and a Failure of Separation on the third Level }. Annals of Pure and Applied Logic 176 (3), 2026.

 
    \bibitem{Uniformization} S. Hoffelner, \textit{Forcing the $\Pi^1_n$-Uniformization Property}, arXiv:2103.11748.
    
    
   \bibitem{Uniformization with wellorder} S. Hoffelner, \textit{A Universe with a $\Delta^1_n$-definable Well-order of the Reals, $\mathsf{CH}$ and $\Pi^1_n$-Uniformization}, arXiv:2506.21778.

    \bibitem{DelfinoSigma} S. Hoffelner and S. M\"uller, \textit{On a local variant of the 12th Delfino problem -- the $\Sigma$-side}, arXiv:2606.13199.

    \bibitem{DelfinoPi} S. Hoffelner, \textit{On a local variant of the 12th Delfino problem -- the $\Pi$-side}, arXiv:2606.13210.
    
    
    \bibitem{Separation without reduction} S. Hoffelner, \textit{A Failure of $\Pi^1_{n+3}$-Reduction in the Presence of $\Sigma^1_{n+3}$-Separation}, arXiv:2312.02540.

     

   

   
   
    \bibitem{Jech} T. Jech \textit{Set Theory. Third Millenium Edition.} Springer 2006.
    
    
    \bibitem{JensenSolovay}
    
    R. Jensen and R. Solovay \textit{Some Applications of Almost Disjoint Sets.}
    
    Studies in Logic and the Foundations of Mathematics
    
    Volume 59, pp. 84-104, 1970.
    
    
     \bibitem{Kechris} 
A. Kechris 
\textit{Classical Descriptive Set Theory}.  
Springer 1995.
    
    
    
   
     
    
     \bibitem{MS} D. Martin and J. Steel \textit{A Proof of Projective Determinacy.} Journal of the American Mathematical Society (2), pp.71-125, 1989.
    
    
    \bibitem{Moschovakis}
    
    Y. Moschovakis \textit{Descriptive Set Theory.} Mathematical Surveys and Monographs 155, AMS.
    
    
    \bibitem{Moschovakis2}
    
    Y. Moschovakis \textit{Uniformization in a playful Universe.}  Bulletin of the  American Mathematical Society 77, no. 5, 731-736, 1971.


   \bibitem{Schindler}
   R. Schindler \textit{Set Theory: Exploring Independence and Truth}. Springer 2014.
    \end{thebibliography}
    \end{document}